\documentclass{article}

\usepackage[final]{neurips_2022}

\usepackage[utf8]{inputenc} 
\usepackage[T1]{fontenc} 
\usepackage{hyperref} 
\usepackage{url} 
\usepackage{booktabs} 
\usepackage{amsfonts} 
\usepackage{nicefrac} 
\usepackage{microtype} 
\usepackage{xcolor} 

\usepackage{graphicx}
\usepackage{algorithm}
\usepackage{algpseudocode}
\usepackage{amsmath}
\usepackage{amssymb}
\usepackage{mathtools}
\usepackage{tabularx}
\usepackage{multirow}
\usepackage{wrapfig}
\usepackage{amsthm}
\usepackage{arydshln}
\usepackage{subcaption}
\usepackage[version=4]{mhchem}

\DeclareMathOperator*{\argmin}{arg\,min}
\DeclareMathOperator*{\minimize}{minimize}

{
	
}

\theoremstyle{definition}
\newtheorem{definition}{Definition}

\newtheorem{theorem}{Theorem}

\newtheorem{lemma}{Lemma}

\usepackage{stackengine}
\def\delequal{\mathrel{\ensurestackMath{\stackon[1pt]{=}{\scriptstyle\Delta}}}}

\usepackage{array}
\newcolumntype{P}[1]{>{\centering\arraybackslash}p{#1}}

\title{Improved Imaging by Invex Regularizers\\ with Global Optima Guarantees}

\author{%
	Samuel Pinilla$^{1,2}$, Tingting Mu$^{3}$, Neil Bourne$^{2}$, Jeyan Thiyagalingam$^{1}$\\
	$^{1\thanks{Rutherford Appleton Laboratory.} \hspace{0.5em}}$Scientific Computing Department, Science and Technology Facilities Council, Harwell, UK\\
	$^{2}$University of Manchester at Harwell, UK \\
	$^{3}$Computer Science, University of Manchester, UK\\
	\{\texttt{samuel.pinilla,t.jeyan}\}\texttt{@stfc.ac.uk}\\
	\{\texttt{tingting.mu,neil.bourne}\}\texttt{@manchester.ac.uk}
}

\begin{document}
	
	\maketitle
	\begin{abstract}
		Image reconstruction enhanced by regularizers, e.g., to enforce sparsity, low rank or smoothness priors on images, has many successful applications in vision tasks such as computer photography, biomedical and spectral imaging. It has been well accepted that non-convex regularizers normally perform better than convex ones in terms of the reconstruction quality. But their convergence analysis is only established to a critical point, rather than the global optima. To mitigate the loss of guarantees for global optima, we propose to apply the concept of \textit{invexity} and provide the first list of proved invex regularizers for improving image reconstruction. Moreover, we establish convergence guarantees to global optima for various advanced image reconstruction techniques after being improved by such invex regularization. To the best of our knowledge, this is the first practical work applying invex regularization to improve imaging with global optima guarantees. To demonstrate the effectiveness of invex regularization, numerical experiments are conducted for various imaging tasks using benchmark datasets. \vspace{-1em}
	\end{abstract}

	\section{Introduction} 
	Image reconstruction (restoration) enhanced by regularizers has a wide application in vision tasks such as computed tomography \cite{sidky2012convex,zhang2022opk_snca}, optical imaging \cite{meinhardt2017learning,afonso2010fast}, magnetic resonance imaging \cite{ehrhardt2016multicontrast,fessler2020optimization}, computer photography \cite{rostami2021power,heide2013high}, biomedical and spectral imaging \cite{wang2019hyperspectral,zhang2013parallel}. In general, an image reconstruction task can be formulated as the solution of the following optimization problem:
	\begin{align}
		\minimize_{\boldsymbol{x}\in \mathbb{R}^{n}} \hspace{0.5em}F(\boldsymbol{x}) = f(\boldsymbol{x}) + g(\boldsymbol{x}).
		\label{eq:basicProblem}
	\end{align}
	Here $f(\boldsymbol{x})$ models a data fidelity term, which usually corresponds to an error loss for image reconstruction, and is assumed to be differentiable. The other function $g(\boldsymbol{x})$ acts as a regularizer which can be non-smooth. It imposes image priors such as sparsity, low rank or smoothness \cite{monga2017handbook}. The use of an appropriate regularizer plays an important role in obtaining robust reconstruction results.
	
	Convex regularization has been popular in the last decade \cite{monga2017handbook,sun2019computed,beck2009fast,liu2016projected,soldevila2016computational}, because it can result in guaranteed global optima. The most well-known examples include the $\ell_{1}$-norm and nuclear norm, which are the continuous and convex surrogates of the $\ell_{0}$-pseudo norm and rank, respectively \cite{fu2014low}. Although convex regularizers have demonstrated their success in signal/image processing, biomedical informatics and computer vision applications \cite{beck2009fast,shevade2003simple,wright2008robust,ye2012sparse}, they are suboptimal in many cases, as they promote sparsity and low rank only under very limited conditions (more measurements from the scene are needed \cite{candes2008enhancing,zhang2010analysis}). To address such limitations, non-convex regularizers have been proposed. For instance, several interpolations between the $\ell_{0}$-pseudonorm and the $\ell_{1}$-norm have been explored including the $\ell_{p}$-quasinorms (where $0<p<1$) \cite{marjanovic2012l_q}, Capped-$\ell_{1}$ penalty \cite{zhang2007surrogate}, Log-Sum Penalty \cite{candes2008enhancing}, Minimax Concave Penalty \cite{zhang2010nearly}, Geman Penalty \cite{geman1995nonlinear}. However, these non-convex regularizers unfortunately come with the price of losing global optima guarantees. 
	
	Image reconstruction methods based on Eq. \eqref{eq:basicProblem} include model-based approaches that directly solve Eq. \eqref{eq:basicProblem} using well-established optimization techniques, e.g., proximal operators and gradient descent rules \cite{beck2017first,jin2021non,sarao2021analytical}, learning-based approaches that train an inference neural network \cite{zhang2021plug,goodfellow2016deep}, as well as hybrid approaches that draw links between iterative signal processing algorithms and the layer-wise neural network architectures \cite{pinilla2022unfolding,monga2021algorithm}. Many of these exploit non-convex assumptions over $f(\boldsymbol{x})$ and/or $g(\boldsymbol{x})$, for which we present a summary of some commonly used or successful ones in Table \ref{tab:literatureComposite}. The table includes algorithms like the iterative reweighted least squares (IRLS) \cite{mohan2012iterative,ochs2015iteratively}, where the regularizer is a composition between the one-dimensional $\ell_{p}$-quasinorm and the trace of a matrix. In \cite{attouch2013convergence,frankel2015splitting}, the objective function $F(\boldsymbol{x})$ is assumed to form a semi-algebraic or tame optimization problem solved by gradient descent algorithms. In \cite{gong2013general}, the regularizer $g(\boldsymbol{x})$ is assumed to be the subtraction of two convex functions, and the general iterative shrinkage and thresholding (GIST) algorithm is proposed to optimize $F(\boldsymbol{x})$. Lastly, \cite{ochs2014ipiano} assumes non-convex $f(\boldsymbol{x})$ but convex $g(\boldsymbol{x})$ and proposes the inertial proximal (iPiano) algorithm for optimization.
	\begin{table}[t!]
		\centering
		\caption{\small Comparison between the assumptions made in this work for $f(\boldsymbol{x})$, and $g(\boldsymbol{x})$ to be optimized in Eq. \eqref{eq:basicProblem} and the most common/successful assumptions in the state-of-the-art.}
		\footnotesize
		\begin{tabular}{m{3.5cm} m{6cm} m{3cm} }
			\hline
			\hline 
			\textbf{Method name} \vspace{0.4em} & \textbf{Assumption} \vspace{0.4em}& \textbf{Global optimizer} \vspace{0.4em}\\
			\hline 
			IRLS \cite{mohan2012iterative,ochs2015iteratively} & special $f$ and $g$ & No \\
			\hline 
			General descent \cite{attouch2013convergence,frankel2015splitting} & Kurdyka-Łojasiewicz & No \\
			\hline 
			GIST \cite{gong2013general} & nonconvex $f$, $g = g_{1} -g_{2}$, $g_{1},g_{2}$ convex & No \\
			\hline 
			iPiano \cite{ochs2014ipiano} & nonconvex $f$, convex $g$ & No \\
			\hline 
			\textbf{Proposed} & convex $f$, invex $g$ & \textbf{Yes} \\
			\hline
			\hline
		\end{tabular}
		\label{tab:literatureComposite}
		\vspace{-1em}
	\end{table}
	
	For algorithms with the convexity assumptions removed, e.g., those in Table \ref{tab:literatureComposite}, their convergence analysis unfortunately can only be established for a critical point. Ideally, we always prefer algorithms that can find the optimal solution for the target problem. One way to mitigate the loss of guarantees for global optima is by revisiting the concept of \textit{invexity} which was first introduced by Hanson \cite{hanson1981sufficiency}, Craven and Glover \cite{craven1985invex} in the 1980s. What makes this class of functions special is that, for any point where the derivative of a function vanishes (stationary point), it is a global minimizer of the function. Convexity is a special case of invexity. Since 1990s, a lot of mathematical implications for invex functions have been developed, but with the lack of practical applications \cite{zualinescu2014critical}. Examples of the few successful works implementing the invexity theory include \cite{barik2021fair,syed2013invexity,chen2016generalized}. To the best of our knowledge, there is no existing work on the application of invex regularization for imaging. 
	
	In this paper, we focus on image reconstruction problems formulated in the form of Eq. \eqref{eq:basicProblem}, where the data fidelity term $f(\boldsymbol{x})$ is based on the $\ell_{2}$-norm and an invex regularizer $g(\boldsymbol{x})$ is used. Most invex theory research lacks clarity on how to benefit practical applications, and this does not encourage the practitioners to exploit the invex property \cite{zualinescu2014critical}. We aim at filling this gap by providing for the first time concrete and useful invex optimization formulations for imaging applications.
	
	Specifically, we make the following contribution:
	\begin{itemize}
		\item Provide the first list of regularizers with proved invexity that fits optimization problems for imaging applications.
		
		\item Establish convergence guarantees to global optima for three types of advanced image reconstruction techniques enhanced by invex reguarlizers.
		
		\item Empirically demonstrate the effectiveness of invex regularization for various imaging tasks. \vspace{-0.5em}
	\end{itemize}

	\section{Preliminaries} 
	\vspace{-0.5em}
	Throughout this paper, we use boldface lowercase and uppercase letters for vectors and matrices, respectively. The $i$-th entry of a vector $\boldsymbol{w}$, is $\boldsymbol{w}[i]$. For vectors, $\lVert \boldsymbol{w}\rVert_p$ is the $\ell_p$-norm. An open ball is defined as $B(\boldsymbol{x};r) = \left \lbrace \boldsymbol{y}\in \mathbb{R}^{n}: \lVert \boldsymbol{y}-\boldsymbol{x} \rVert_{2}<r \right\rbrace$. The operation $\text{conv}(\mathcal{A})$ represents the convex hull of the set $\mathcal{A}$, and the operation $\text{sign}(w)$ returns the sign of $w$. We use $\sigma_{i}(\boldsymbol{W})$ to denote the $i$-th singular value of $\boldsymbol{W}$ assumed in descending order. 
	
	We present several concepts needed for the development of this paper starting with the definition of a locally Lipschitz continuous function. 
	\begin{definition}[\textbf{Locally Lipschitz Continuity}]
		A function $f:\mathbb{R}^{n}\rightarrow \mathbb{R}$ is locally Lipschitz continuous at a point $\boldsymbol{x}\in \mathbb{R}^{n}$ if there exist scalars $K>0$ and $\epsilon>0$ such that 
		\begin{align}
			\lvert f(\boldsymbol{y}) - f(\boldsymbol{z}) \rvert \leq K\lVert \boldsymbol{y}-\boldsymbol{z} \rVert_{2},
		\end{align}
		for all $\boldsymbol{y},\boldsymbol{z}\in B(\boldsymbol{x},\epsilon)$.
		\label{def:lipschitz}
	\end{definition}
	
	Since the ordinary directional derivative being the most important tool in optimization does not necessarily exist for locally Lipschitz continuous functions, it is required to introduce the concept of subdifferential \cite{B2014} which is calculated in practice as follows. 
	
	\begin{theorem}[\textbf{Subdifferential}]{\cite[Theorem 3.9]{B2014}}
		\label{theo:auxDerivative}
		Let $f:\mathbb{R}^{n}\rightarrow \mathbb{R}$ be a locally Lipschitz continuous function at $\boldsymbol{x}\in \mathbb{R}^{n}$, and define $\Omega_{f} = \{\boldsymbol{x}\in \mathbb{R}^{n}| \text{ f is not differentiable at the point } \boldsymbol{x}\}$. Then the subdifferential of $f$ is given by
		\begin{align}
			\partial f(\boldsymbol{x}) = \text{ conv }&\left( \left \lbrace \boldsymbol{\zeta}\in \mathbb{R}^{n}| \text{ exists } (\boldsymbol{x}_{i})\in \mathbb{R}^{n}\setminus \Omega_{f} \text{ such that } \boldsymbol{x}_{i}\rightarrow \boldsymbol{x} \text{ and }\nabla f(\boldsymbol{x}_{i})\rightarrow \boldsymbol{\zeta} \right \rbrace\right).
		\end{align}
	\end{theorem}
	
	The notion of subdifferential is given for locally Lipschitz continuous functions because it is always nonempty \cite[Theorem 3.3]{B2014}. Based on these, the concept of invex function is presented as follows. 
	\begin{definition}[\textbf{Invexity}]
		\label{def:invex}
		Let $f:\mathbb{R}^{n}\rightarrow \mathbb{R}$ be locally Lipschitz; then $f$ is invex if there exists a function $\eta:\mathbb{R}^{n}\times \mathbb{R}^{n} \rightarrow \mathbb{R}^{n}$ such that
		\begin{align}
			f(\boldsymbol{x})-f(\boldsymbol{y}) \geq \boldsymbol{\zeta}^{T}\eta(\boldsymbol{x},\boldsymbol{y}),
			\label{eq:basicInvex}
		\end{align}
		$\forall \boldsymbol{x},\boldsymbol{y} \in \mathbb{R}^{n}$, $\forall \boldsymbol{\zeta} \in \partial f(\boldsymbol{y})$.
	\end{definition}
	It is well known that a convex function simply satisfies this definition for $\eta(\boldsymbol{x},\boldsymbol{y}) = \boldsymbol{x}-\boldsymbol{y}$. 
	
	The following classical theorem \cite[Theorem 4.33]{mishra2008invexity} makes connection between an invex function and its well-known optimum property that supports the motivation of designing invex regularizers.
	\begin{theorem}[\textbf{Invex Optimality}]{\cite[Theorem 4.33]{mishra2008invexity})}
		Let $f:\mathbb{R}^{n}\rightarrow \mathbb{R}$ be locally Lipschitz. Then the following statements are equivalent.
		\begin{enumerate}
			\item $f$ is invex.
			\item Every point $\boldsymbol{y}\in \mathbb{R}^{n}$ that satisfies $\boldsymbol{0}\in \partial f(\boldsymbol{y})$ is a global minimizer of $f$.
			\item Definition \ref{def:invex} is satisfied for $\eta:\mathbb{R}^{n}\times \mathbb{R}^{n} \rightarrow \mathbb{R}^{n}$ given by
			\begin{align}
				\eta(\boldsymbol{x},\boldsymbol{y})= \left \lbrace \begin{array}{ll}
					\boldsymbol{0} & f(\boldsymbol{x})\geq f(\boldsymbol{y}),\\
					\frac{f(\boldsymbol{x})-f(\boldsymbol{y})}{\lVert \boldsymbol{\zeta}^{*}_{\boldsymbol{y}}\rVert_{2}^{2}}\boldsymbol{\zeta^{*}_{y}} & \text{ otherwise, }
				\end{array}\right.
			\end{align}
			where $\boldsymbol{\zeta^{*}_{y}}$ is an element in $\partial f(\boldsymbol{y})$ of minimum norm.
		\end{enumerate}
		\label{theo:optimal_v0}
	\end{theorem} 
	
	\section{Invex Functions}
	\label{sec:InvexRegul}\vspace{-0.8em}
	We start this section by firstly presenting five examples of invex functions that are useful for imaging applications. Four of these have been labelled as non-convex in existing works \cite{wu2019improved,wen2018survey}. This is the first time that they are formally proved to be invex functions. We prove their invexity by showing they satisfy Statement 2 of Theorem \ref{theo:optimal_v0} (see proof in Appendix \ref{app:invexProof} of supplementary material).
	\begin{lemma} [\textbf{Invex Functions}]
		All of the following functions are invex:
		\begin{align}
			\label{fun1}
			g(\boldsymbol{x}) =& \sum_{i=1}^{n}\left(\lvert\boldsymbol{x}[i] \rvert + \epsilon \right)^{p}, \textmd{for }p\in (0,1) \textmd{ and } \epsilon\geq \left(p(1-p)\right)^{\frac{1}{2-p}}, \\
			\label{fun2}
			g(\boldsymbol{x}) = &\sum_{i=1}^{n}\log(1 + \lvert \boldsymbol{x}[i] \rvert),\\
			\label{fun3}
			g(\boldsymbol{x}) = & \sum_{i=1}^{n}\frac{\lvert \boldsymbol{x}[i]\rvert}{2 + 2\lvert \boldsymbol{x}[i] \rvert},\\
			\label{fun4}
			g(\boldsymbol{x}) = &\sum_{i=1}^{n}\frac{\boldsymbol{x}^{2}[i]}{1 + \boldsymbol{x}^{2}[i]},
		\end{align} 
		\begin{align}
			\label{fun5}
			g(\boldsymbol{x})= &\sum_{i=1}^{n} \log(1+\lvert \boldsymbol{x}[i] \rvert) - \frac{\lvert \boldsymbol{x}[i] \rvert}{2 + 2\lvert \boldsymbol{x}[i] \rvert}.
		\end{align}
		\label{theo:invexProof}
		\vspace{-1em}
	\end{lemma}
	
	We provide further insights of these functions in Section \ref{app:discussionRegu} of Supplemental material. Table \ref{tab:list} summarizes their applications. Specifically, Eq. \eqref{fun1} is known as quasinorm, and has attracted a lot of attention because it has resulted in theoretical improvements for matrix completion and compressive sensing \cite{marjanovic2012l_q,wu2013improved}. The analysis on the quasinorms is valid with and without the constant $\epsilon$. We prefer to add $\epsilon$ in order to formally satisfy the Lipschitz continuity in Definition \ref{def:lipschitz}. Eqs. \eqref{fun2} and \eqref{fun3} enhance the convex $\ell_{1}$-norm regularizer, and they have significantly improved image denoising \cite{wu2019improved}. Eq. \eqref{fun4} has been used as the loss function to improve support vector classification \cite{zhuang2019surrogate}. 
	
	\begin{table}[ht]
		\centering
		\caption{List of invex functions studied in this work.}
		\resizebox{0.8\textwidth}{!}{
			\begin{tabular}{ m{2.5cm} m{2.5cm} m{5cm} }
				\hline
				\hline 
				\textbf{Reference} \vspace{0.4em} & \textbf{Invex function} \vspace{0.4em} & \textbf{Application} \vspace{0.4em}\\
				\hline
				\cite{marjanovic2012l_q,mohan2012iterative,zhang2018sparse,lin2019efficient} & Eq. \eqref{fun1} & Matrix completion \\
				\hline
				\cite{candes2008enhancing,gong2013general,yao2016efficient,wen2018survey} & Eq. \eqref{fun2} & Enhancing compressive sensing \\
				\hline 
				\cite{wu2019improved,hu2021low,wang2022accelerated} & Eq. \eqref{fun3} & Image denoising \\
				\hline
				\cite{zhuang2019surrogate} & Eq. \eqref{fun4} & Support vector classification \\
				\hline
				Proposed & 
				Eq. \eqref{fun5} & Compressive sensing \\
				\hline
		\end{tabular}}
		\label{tab:list}
	\end{table} 
	
	We propose the last function in Eq. \eqref{fun5} by the subtraction between Eq. \eqref{fun2} and Eq. \eqref{fun3}. This design is motivated by the optimization framework in \cite{gong2013general} where the regularization term is assumed to be the subtraction of two convex functions (see GIST in Table \ref{tab:literatureComposite}). This has been found to be highly successful in imaging applications (see the survey \cite{wen2018survey}). But until now there is no evidence that this subtraction produces another convex function (if exists) potentially useful in imaging applications. Therefore, we propose this example to show that at least this is possible in the invex case. 	
	
	Additionally, we present another way of constructing an invex function in the following lemma. It establishes that an invex function $f:\mathbb{R}^{m}\rightarrow \mathbb{R}$ composed with an affine mapping $\boldsymbol{H}\boldsymbol{x}-\boldsymbol{b}$ for $\boldsymbol{H}\in \mathbb{R}^{m\times n}$, $\boldsymbol{x}\in \mathbb{R}^{n}$ and $\boldsymbol{b}\in \mathbb{R}^{m}$, is also invex if $\boldsymbol{H}$ is full row-rank. This condition on $\boldsymbol{H}$ is a mild assumption, because we show in Section \ref{sec:invexImag} imaging application examples that satisfy this criterium.
	\begin{lemma}[\textbf{Affine Invex Construction}]
		Let $f:\mathbb{R}^{m}\rightarrow \mathbb{R}$ be a continuously differentiable invex function, $\boldsymbol{H}\in \mathbb{R}^{m\times n}$ have full row rank, and $\boldsymbol{b}\in \mathbb{R}^{m}$ be a vector. Then the function $h(\boldsymbol{x}) = f(\boldsymbol{H}\boldsymbol{x}-\boldsymbol{b})$ is invex.
		\label{theo:invexComposited}
	\end{lemma}
	
	Similar to Lemma \ref{theo:invexProof}, it is proved by showing that the composed function satisfies Statement 2 of Theorem \ref{theo:optimal_v0} (see Appendix \ref{app:invexComposited} of supplementary material). Eq. \eqref{fun4} is an example of such an invex construction that satisfies the continuously differentiable assumption in Lemma \ref{theo:invexComposited}, easily verified its proof in Appendix \ref{app:invexComposited}. A practical implication of Lemma \ref{theo:invexComposited} for imaging applications appears when we want to solve linear system of equations (e.g. \cite{zhuang2019surrogate}). We demonstrate an application of this kind of invex construction in Section \ref{sub:PnP} to improve a widely used image restoration framework.
	
	\section{Invex Imaging Examples, Algorithms and Convergence Analysis}
	\label{sec:invexImag} 
	In this section, we demonstrate the use of invex regularizers to improve some advanced imaging methodologies. To benefit both practitioners and theory development, we present practical invex imaging algorithms and prove their convergence guarantees to global optima which was only possible for convex functions.	
	
	\subsection{Image Denoising}
	\label{sub:denoising}
	Image denoising plays a critical role in modern signal processing systems since images are inevitably contaminated by noise during acquisition, compression, and transmission, leading to distortion and loss of image information \cite{fan2019brief}. Plenty of denoising methods exist, originating from a wide range of disciplines such as probability theory, statistics, partial differential equations, linear and nonlinear filtering, spectral and multiresolution analysis, also classical machine learning and deep learning \cite{mahmoudi2005fast,krull2019noise2void,fan2019brief}. All these methods rely on some explicit or implicit assumptions about the true (noise-free) signal in order to separate it properly from the random noise. 
	
	One of the most successful assumptions is that a signal can be well approximated by a linear combination of few basis elements in a transform domain \cite{dabov2007image,elad2006image}. Under this assumption, a denoising method can be implemented as a two-step procedure: i) to obtain high-magnitude transform coefficients that convey mostly the true-signal energy, ii) to discard the transform coefficients which are mainly due to noise. Typical choices for the first step are the wavelet, cosine transforms, and principal component analysis (PCA) \cite{dabov2007image,elad2006image,cai2014data}. The second step is seen as a filtering procedure that is formally modelled as a proximal optimization problem~\cite{parikh2014proximal}
	\begin{align}
		\text{Prox}_{g}(\boldsymbol{u}) = \argmin_{\boldsymbol{x}\in \mathbb{R}^{n}} \left( g(\boldsymbol{x}) + \frac{1}{2}\lVert \boldsymbol{x} -\boldsymbol{u}\rVert_{2}^{2}\right),
		\label{eq:prox1}
	\end{align}
	where $g(\boldsymbol{x})$ acts as a regularization term, and $\boldsymbol{u}$ represents the noisy transform coefficients. In fact, the usefulness of Eq. \eqref{eq:prox1} is not just limited to denoising, but other imaging problems like computer tomography \cite{jorgensen2021core}, optical imaging \cite{sun2019regularized}, biomedical and spectral imaging \cite{sun2019online}. In general, global optima guarantees in Eq. \eqref{eq:prox1} is restricted to convex $g(\boldsymbol{x})$, e.g., $\ell_{1}$-norm. 
	
	We improve this important proximal operator by incorporating invex regularizers. Specifically, using those invex functions $g(\boldsymbol{x}) $ as listed in Table \ref{tab:list}, global minimization is achieved in Eq. \eqref{eq:prox1}. The result is presented in the following theorem: 	
	\begin{theorem}[\textbf{Invex Proximal}]
		Consider the optimization problem in Eq. \eqref{eq:prox1} for all functions in Table \ref{tab:list}. Then the following holds:
		\begin{enumerate}
			\item The function $h(\boldsymbol{x}) = g(\boldsymbol{x}) + \frac{1}{2}\lVert \boldsymbol{x} -\boldsymbol{u}\rVert_{2}^{2}$ is convex (therefore invex).
			\item The resolvent operator of the proximal is $(\mathbf{I} + \partial g)^{-1}$ and it is treated as a singleton because it always maps to a global optimizer.
		\end{enumerate}
		\label{theo:proximalProof}
	\end{theorem}
	It is classically known that the sum of two invex functions is not necessarily invex in general \cite{mishra2008invexity}. Therefore, presenting examples like above, where the sum of $f(\boldsymbol{x}) $ and $g(\boldsymbol{x}) $ is invex, is important to both invexity and imaging communities. We present the proof of Theorem \ref{theo:proximalProof} and provide the solution to Eq. \eqref{eq:prox1} for each function in Table \ref{tab:list} in Appendix \ref{app:proximalProof} of supplementary material. 
	
	\subsection{Image Compressive Sensing} 
	\label{sub:CSAlg}
	Image \textit{compressive sensing} has been extensively exploited in areas such as microscopy, holography, optical imaging and spectroscopy \cite{arce2013compressive,jerez2020fast,guerrero2020phase}. It is an inverse problem that aims at recovering an image $\boldsymbol{f}\in \mathbb{R}^{n}$ from its measurement data vector $\boldsymbol{b} = \boldsymbol{\Phi}\boldsymbol{f}$, where $\boldsymbol{\Phi} \in \mathbb{R}^{m\times n}$ is the image acquisition matrix ($m<n$). Since $m<n$, compressive sensing assumes $\boldsymbol{f}$ has a $k$-sparse representation $\boldsymbol{x}\in \mathbb{R}^{n}$ ($k\ll n$ non-zero elements) in a basis $\boldsymbol{\Psi}\in \mathbb{R}^{n\times n}$, that is $\boldsymbol{f}=\boldsymbol{\Psi}\boldsymbol{x}$, in order to ensure uniqueness under some conditions. Examples of this sparse basis $\boldsymbol{\Psi}$ in imaging are the Wavelet (also Haar Wavelet) transform, cosine and Fourier representations \cite{foucart13}. Hence, one can work with the abstract model $\boldsymbol{b}= \boldsymbol{\Phi}\boldsymbol{\Psi}\boldsymbol{x}=\boldsymbol{H}\boldsymbol{x}$, where $\boldsymbol{H}$ encapsulates the product between $\boldsymbol{\Phi}$, and $\boldsymbol{\Psi}$, with $\ell_{2}$-normalized columns \cite{arce2013compressive,candes2008introduction}. Under this setup, compressive sensing enables to recover $\boldsymbol{x}$ using much lesser samples than what are predicted by the Nyquist criterion \cite{candes2008introduction}. The task formulation is
	\begin{align}
		\minimize_{\boldsymbol{x}\in \mathbb{R}^{n}} \hspace{0.5em} & f(\boldsymbol{x}) + \lambda g(\boldsymbol{x}) = \frac{1}{2}\lVert \boldsymbol{H}\boldsymbol{x} - \boldsymbol{b} \rVert_{2}^{2} + \lambda g(\boldsymbol{x}),
		\label{eq:problem4}
	\end{align}
	where $\lambda\in (0,1]$ is a typical choice in practice. When the regularizer $g(\boldsymbol{x})$ takes the convex form of $\ell_{1}$-norm, and when the sampling matrix $\boldsymbol{H}$ satisfies the \textit{restricted isometry property} (RIP) for any $k$-sparse vector $\boldsymbol{x}\in \mathbb{R}^{n}$, i.e., $(1-\delta_{2k})\lVert \boldsymbol{x} \rVert_{2}^{2} \leq \lVert \boldsymbol{H}\boldsymbol{x} \rVert_{2}^{2} \leq (1+\delta_{2k})\lVert \boldsymbol{x} \rVert_{2}^{2}$ for $\delta_{2k}<\frac{1}{3}$ \cite[Theorem 6.9]{foucart13}, it has been proved that $\boldsymbol{x}$ can be exactly recovered by solving Eq. \eqref{eq:problem4} \cite{candes2006robust}.
	
	We are interested in invex regularizers. It has been proved that, when $g(\boldsymbol{x})$ takes the particular invex form in Eq. \eqref{fun1}, $\boldsymbol{x}$ can be exactly recovered by solving Eq. \eqref{eq:problem4} \cite{wu2013improved}. Below we further generalize this result to all the invex functions as listed in Table \ref{tab:list}. The generalized result is presented in Theorem \ref{theo:ourCS}.
	\begin{theorem}[\textbf{Invex Image Compressive Sensing}]
		Assume $\boldsymbol{H}\boldsymbol{x}=\boldsymbol{b}$, where $\boldsymbol{x}\in \mathbb{R}^{n}$ is $k$-sparse, the matrix $\boldsymbol{H}\in \mathbb{R}^{m\times n}$ ($m<n$) with $\ell_{2}$-normalized columns that satisfies the RIP condition for any $k$-sparse vector, and $\boldsymbol{b}\in \mathbb{R}^{m}$ is a noiseless measurement vector. If $g(\boldsymbol{x})$ in Eq. \eqref{eq:problem4} takes the form of the functions in Table \ref{tab:list}, then the following holds:
		\begin{enumerate}
			\item The objective function $\frac{1}{2}\lVert \boldsymbol{H}\boldsymbol{x} - \boldsymbol{b} \rVert_{2}^{2} + \lambda g(\boldsymbol{x})$ is invex.
			\item $\boldsymbol{x}$ can be exactly recovered by solving Eq. \eqref{eq:problem4} i.e. only global optimizers exist. When $g(\boldsymbol{x})$ takes the form of Eq. \eqref{fun4}, extra mild conditions on $\boldsymbol{x}$ are needed.
		\end{enumerate}
		\label{theo:ourCS}
	\end{theorem}
	We clarify that if $\boldsymbol{H}$ satisfies the mentioned RIP, then each sub-matrix with $k$-columns of $\boldsymbol{H}$, selected according to indices of the nonzero elements of the $k$-sparse signal is a full row-rank matrix. This result is important to invex community, because it supports the validity of Lemma \ref{theo:invexComposited} to build invex functions with affine mappings. Additionally, we present another proved form of function sum that can result in an invex function, i.e., the sum of $g(\boldsymbol{x})$ and the $\ell_{2}$-norm composed with the affine mapping $\boldsymbol{H}\boldsymbol{x}-\boldsymbol{b}$. The complete proof is provided in Appendix \ref{app:ourCS} of supplementary material. 
	
	Next, we present different algorithms to solve Eq. \eqref{eq:problem4} using invex $g(\boldsymbol{x})$ as in Table \ref{tab:list}. We select a few of the most important and successful image reconstruction techniques to start from, and develop their invex extensions. Taking advantage of the invex property, we prove convergence to global minimizers for each extended algorithm, which is unexplored up to date.\vspace{-0.5em}
	
	\begin{algorithm}[ht]
		\caption{Accelerated Proximal Gradient}
		\label{alg:invexProximal}
		\begin{algorithmic}[1]
			\State{\textbf{input}: Tolerance constant $\epsilon\in (0,1)$, initial point $\boldsymbol{x}^{(0)}$, and number of iterations $T$.}
			\State{\textbf{initialize}: $\boldsymbol{x}^{(1)}=\boldsymbol{x}^{(0)}=\boldsymbol{z}^{(0)}, r_{1}=1,r_{0}=0, \alpha_{1},\alpha_{2}< \frac{1}{L}$, and $\lambda \in (0,1]$}
			\For{$t=1$ to $T$}
			\State{$\boldsymbol{y}^{(t)}= \boldsymbol{x}^{(t)} + \frac{r_{t-1}}{r_{t}}(\boldsymbol{z}^{(t)}-\boldsymbol{x}^{(t)}) + \frac{r_{t-1}-1}{r_{t}}(\boldsymbol{x}^{(t)}- \boldsymbol{x}^{(t-1)})$}
			\State{$\boldsymbol{z}^{(t+1)}=\text{prox}_{\alpha_{2} \lambda g}(\boldsymbol{y}^{(t)} - \alpha_{2}\nabla f(\boldsymbol{y}^{(t)}))$}
			\State{$\boldsymbol{v}^{(t+1)}=\text{prox}_{\alpha_{1} \lambda g}(\boldsymbol{x}^{(t)} - \alpha_{1}\nabla f(\boldsymbol{x}^{(t)}))$}
			\State{$r_{t+1}=\frac{\sqrt{4(r_{t})^{2}+1}+1}{2}$}
			\State{$\boldsymbol{x}^{(t+1)}=\left \lbrace\begin{array}{ll}
					\boldsymbol{z}^{(t+1)}, & \text{ if }f(\boldsymbol{z}^{(t+1)} ) + \lambda g(\boldsymbol{z}^{(t+1)} )\leq f(\boldsymbol{v}^{(t+1)} ) + \lambda g(\boldsymbol{v}^{(t+1)} )\\
					\boldsymbol{v}^{(t+1)}, & \text{ otherwise }
				\end{array}\right.$}
			\EndFor
			\State{\textbf{return:} $\boldsymbol{x}^{(T)}$}
		\end{algorithmic}
	\end{algorithm}\vspace{-0.5em}
	
	\subsubsection{Accelerated Proximal Gradient Algorithm} 
	The accelerated proximal gradient (APG) method \cite{li2015accelerated} has been shown to be effective solving Eq. \eqref{eq:problem4}, achieving better imaging quality in less iterations than its predecessors \cite{beck2009fast,frankel2015splitting,gong2013general,ochs2014ipiano,boct2016inertial}, and been frequently used by recent imaging works \cite{wang2022accelerated,mai2022energy,zhang2022continual,ge2022fast}. Its convergence to global optima is only guaranteed for convex loss \cite{li2015accelerated}. For non-convex cases, convergence to a critical point has been stated \cite{li2015accelerated}. Its pseudo-code for solving Eq. \eqref{eq:problem4} is provided in Algorithm \ref{alg:invexProximal}. 
	
	Taking advantage that the loss function $f(\boldsymbol{x}) + \lambda g(\boldsymbol{x})$ in Eq. \eqref{eq:problem4} is invex, and the uniqueness result in Theorem \ref{theo:proximalProof}, we formally extend APG in the following lemma stating that the sequence $\left\{\boldsymbol{x}^{(t+1)}\right\}$ generated by Algorithm \ref{alg:invexProximal} converges to a global minimizer of Eq. \eqref{eq:problem4}. 
	\begin{lemma}[\textbf{Invex APG}]
		\label{lem:convergeAPG}
		Under the setup of Theorem \ref{theo:ourCS} and using $L=\sigma_{1}\left(\boldsymbol{H}^{T}\boldsymbol{H}\right)$ (maximum singular value), the sequence $\left\{\boldsymbol{x}^{(t)}\right\}_{t=0}^{T-1}$ generated by Algorithm \ref{alg:invexProximal} converges to a global minimizer.
	\end{lemma} 
	To prove Lemma \ref{lem:convergeAPG}, we apply the Statement 2 of Theorem \ref{theo:optimal_v0} to Eq. \eqref{eq:problem4} and the unicity of the proximal operators for functions in Table \ref{tab:list}. The proof is provided in Appendix \ref{app:lemAPG} of supplementary material.

	\subsubsection{Plug-and-play with Deep Denoiser Prior} 
	\label{sub:PnP}
	Plug-and-play (PnP) is a powerful framework for regularizing imaging inverse problems \cite{sun2019online} and has gained popularity in a range of applications in the context of imaging inverse problems \cite{zhang2021plug,sun2019online,wei2022tfpnp,kamilov2022plug,hu2022monotonically}. It replaces the proximal operator in an iterative algorithm with an image denoiser, which does not necessarily have a corresponding regularization objective. This implies that the effectiveness of PnP goes beyond standard proximal algorithms such as primal-dual splitting \cite{ono2017primal,kamilov2017plug,zha2022simultaneous}. It has guarantees to a fixed point only when convex objective functions are employed \cite{kamilov2017plug}.
	
	To apply the PnP framework, we modify Algorithm \ref{alg:invexProximal} by replacing the proximal operator (Line 6 in its pseudo-code) with a neural network based denoiser Noise2Void \cite{krull2019noise2void}, resulting in 
	\begin{align}
		\boldsymbol{v}^{(t+1)}=\text{Noise2Void}\left(\boldsymbol{x}^{(t)} - \alpha_{1}\nabla f\left(\boldsymbol{x}^{(t)}\right)\right).
		\label{eq:PnPVariant}
	\end{align}
	
	The complete pseudo-code is presented in Algorithm \ref{alg:invexPnP} of Appendix \ref{app:PnP} in supplemental material. We remark that in Algorithm \ref{alg:invexPnP}, Line 5 of the  Algorithm \ref{alg:invexProximal} is retained to allows the comparison between regularizers (invex and convex). More specifically, Line~5 computes the proximal step, while Line~6 relies on a neural network for the same purpose \eqref{eq:PnPVariant}. This  offers an avenue for simultaneously exploiting both the model-based and data-driven approaches. The output of Algorithm 3 is a close estimation to the solution of Eq. \eqref{eq:problem4} \cite{kamilov2017plug}. The benefit of using this denoiser is that it does not require clean target images in order to be trained. We present the following convergence result for this modified algorithm under the assumption of $f(\boldsymbol{x})$ in Eq. \eqref{eq:problem4} being invex which is a generalization of \cite{kamilov2017plug} (restricted to convex functions only).
	\begin{lemma}[\textbf{Invex Plug-and-play}]
		Assume $f(\boldsymbol{x})$ in Eq. (\ref{eq:problem4}) is invex with Lipschitz continuous gradient, and a denoiser $d:\mathbb{R}^{n}\rightarrow \mathbb{R}$. Under the setup of Theorem \ref{theo:ourCS} and some mild conditions on $d$, the sequence $\left\{\boldsymbol{x}^{(t)}\right\}_{t=0}^{T}$ generated by Algorithm \ref{alg:invexProximal} satisfies
		\begin{align}
			\frac{1}{T}\sum_{t=1}^{T}\left \| \boldsymbol{x}^{(t)} - d\left(\boldsymbol{x}^{(t)}-\alpha_{1}\nabla f\left(\boldsymbol{x}^{(t)}\right)\right) \right\|_{2}^{2} \leq \frac{2}{T}\left(\frac{1+\kappa}{1-\kappa}\right) \left\| \boldsymbol{x}^{(0)} - \boldsymbol{x}^{*} \right\|_{2}^{2},
			\label{eq:fix}
		\end{align}
		for any $\boldsymbol{x}^{*}=d(\boldsymbol{x}^{*}- \alpha_{1}\nabla f(\boldsymbol{x}^{*}))$ (fixed point) and for some $\kappa\in (0,1)$.
		\label{theo:PnP}
	\end{lemma}
	Eq. \eqref{eq:fix} guarantees that $\left\{\boldsymbol{x}^{(t)}\right\}_{t=0}^{T}$ is arbitrarily close to the set of fixed points of $d(\cdot)$, which is considered a close estimation to the solution of Eq. \eqref{eq:problem4} \cite{kamilov2017plug}. Its proof is provided in Appendix \ref{app:PnP} of supplementary material. Eq. \eqref{fun4} is an example satisfying assumption required in Lemma \ref{theo:PnP}.
	
	\subsubsection{Unrolling} 
	\label{unrolling}
	The \textit{unrolling} or \textit{unfolding} framework is another imaging strategy for solving Eq. \eqref{eq:problem4}. It offers a systematic connection between iterative algorithms used in signal processing and the neural networks \cite{pinilla2022unfolding,monga2021algorithm,hu2020iterative}. Unrolled neural networks become popular due to their potential in developing efficient and high-performing network architectures from reasonably sized training sets \cite{chowdhury2021unfolding,naimipour2020upr}. A folded version of the proximal gradient algorithm is presented in Algorithm \ref{alg:unrolling}. Particularly, existing works \cite{chen2018theoretical,liu2019alista} have shown that the efficiency of Algorithm \ref{alg:unrolling} can be improved by simulating a recurrent neural network so that its layers mimic the iterations in Line 4 of Algorithm \ref{alg:unrolling}. Specifically, each $\boldsymbol{x}^{(t+1)}$ constitutes one linear operation which models a layer of the network, followed by a proximal operation that models the activation function. Thus, one forms a deep network by mapping each iteration to a network layer and stacking the layers together to learn $\boldsymbol{H}, \alpha_{t}$, and $\boldsymbol{x}^{(t)}$ for all $t$ which is equivalent to executing an iteration of Algorithm \ref{alg:unrolling} multiple times. Their study was conducted only for $g(\boldsymbol{x})$ in the form of $\ell_{1}$-norm. 
	
	\begin{algorithm}[ht]
		\caption{Folded Proximal Gradient Algorithm}
		\label{alg:unrolling}
		\begin{algorithmic}[1]
			\State{\textbf{input}: initial point $\boldsymbol{x}^{(0)}$, number of iterations $T$}
			\State{\textbf{initialize}: $\alpha_{t}< \frac{2}{L+2}$, and $\lambda \in (0,1]$}
			\For{$t=0$ to $T$}
			\State{$\boldsymbol{x}^{(t+1)}=\text{prox}_{\alpha_{t} \lambda g}(\boldsymbol{x}^{(t)} - \alpha_{t}\boldsymbol{H}^{T}(\boldsymbol{H}\boldsymbol{x}^{(t)}-\boldsymbol{b}))$}
			\EndFor
			\State{\textbf{return:} $\boldsymbol{x}^{(T)}$}
		\end{algorithmic}
	\end{algorithm}
	
	Convergence guarantees to global optima for Algorithm \ref{alg:unrolling} has been established in \cite{beck2009fast}, but it is restricted to convex objective functions. Therefore, due to the success and importance of unrolling we aim to extend the global optima guarantees of Algorithm \ref{alg:unrolling} to invex objectives, and present the results in the following lemma:
	\begin{lemma}[\textbf{Invex Unrolling}]
		\label{lem:convergeUnrolling}
		Under the setup of Theorem \ref{theo:ourCS} and using $L=\sigma_{1}\left(\boldsymbol{H}^{T}\boldsymbol{H}\right)$ (maximum singular value) and $\alpha_{t}<\frac{2}{L+2}$, the sequence $\left\{\boldsymbol{x}^{(t)}\right\}_{t=0}^{T-1}$ generated by Algorithm \ref{alg:unrolling} converges to a global minimizer.
	\end{lemma}
	The key to proving Lemma \ref{lem:convergeUnrolling} relies on the uniqueness result of the proximal operator for functions in Table \ref{tab:list} as stated in Theorem \ref{theo:proximalProof}. The proof is presented in Appendix \ref{app:unrolling} of supplementary material. Such results confirm that the invex unrolled network of Algorithm \ref{alg:unrolling}, which uses the proximal operators of invex mappings as the activation functions, can reach the optimal solution during training. 
	
	\section{Experiments and Results}
	\label{others}
	A number of datasets have been merged to formulate one unique dataset for our training and evaluation purposes. These are DIV2K super-resolution~\cite{agustsson2017ntire}, the McMaster~\cite{zhang2011color}, Kodak~\cite{kodak}, Berkeley Segmentation (BSDS 500) \cite{MartinFTM01}, Tampere Images (TID2013) \cite{ponomarenko2013color} and the Color BSD68 \cite{martin2001database} datasets. We conduct various experiments to study the performance of those invex regularizers as listed in Table \ref{tab:list} in non-ideal conditions. We compare them against the state-of-the-art methods originally developed for convex regularizers ($\ell_{1}$-norm) ensuring global optima. When neural network training is involved, we take a total of 900 images which are randomly divided into a training set of 800 images, a validation set of 55 images, and a test set of 45 images. For all the experiments, the images are scaled into the range $[0,1]$. For the invex regularizer in Eq. \eqref{fun1}, we vary the value of $p$.
	
	\subsection{Image Compressive Sensing Experiments} 	
	We assess signal reconstruction, in these experiments, by averaging the peak-signal-to-noise-ratio (PSNR) in dB over the testing image set. We consider additive white Gaussian noise in the measurements data vector with three different levels of SNR (Signal-to-Noise Ratio) = 20, 30, and $\infty$ (noiseless case). For Algorithm \ref{alg:invexProximal} and its plug-and-play variant, the parameters $\lambda,\alpha_{1}$, and $\alpha_{2}$ were chosen to be the best for each analyzed function determined by cross validation, and the initial point $\boldsymbol{x}^{(0)}$ was the blurred image $\boldsymbol{b}$. The results are summarized in Table \ref{tab:globalResults}, where the best and least efficient among invex functions is highlighted in boldface and underscore, respectively. Additional results are reported in Appendix \ref{app:newResults} of supplemental material for each experiment, using the structural similarity index measure to assess imaging quality.
	
	
	\textbf{Experiment 1} studies the effect of different invex regularizers, the Smoothly Clipped Absolute Deviation (SCAD) \cite{fan2001variable}, and the Minimax Concave Penalty (MCP) \cite{zhang2010nearly}, under Algorithm \ref{alg:invexProximal}. A deconvolution problem is studied to formulate Eq. \eqref{eq:problem4} which is an important problem in signal processing due to imperfect artefacts in physical setups such as mismatch, calibration errors, and loss of contrast \cite{yeh2015experimental}. To compare, the used state-of-the-art methods that employ convex regularization are the Total Variation Minimization by Augmented Lagrangian (TVAL3) \cite{li2013efficient}, and the fast iterative shrinkage-thresholding algorithm (FISTA) \cite{beck2009fast} which ensures global optima. Further, to comparing with convolutional neural networks methodologies, the non-iterative reconstruction methodology ReconNet \cite{kulkarni2016reconnet} is used. To model this problem, all pixels of the testing set are fixed to $256\times 256$ pixels. The images went through a Gaussian blur of size $9\times 9$ and standard deviation $4$, followed by an additive zero-mean white Gaussian noise. The sensing matrix $\boldsymbol{H}$ is built as $\boldsymbol{H}=\boldsymbol{\Phi}\boldsymbol{\Psi}$ (for all methods except ReconNet), where $\boldsymbol{\Phi}$ represents the blur operator over the images and $\boldsymbol{\Psi}$ is the inverse of a three stage Haar wavelet transform. This experiment is extremely ill-conditioned, where the condition number of $\boldsymbol{H}^{T}\boldsymbol{H}$ is significantly higher than 1. This means that in practice the RIP condition is not guaranteed. To achieve a fair comparison, the number of iterations was fixed for all functions as $T=800$. The deconvolution problem follows a compressive sensing setup because the Gaussian filter remove high frequency information of the input image.
	
	In the case of ReconNet, we follow existing setting in \cite{kulkarni2016reconnet}. For the learning of ReconNet, we extract patches of size $33\times 33$ from the noisy blurred training image set, and we train it using the Adam optimization algorithm and a learning rate $5\times 10^{-4}$ for 512 epochs with a batch size of $128$.
	
	\textbf{Experiment 2} studies the invex regularizers under the plug-and-play modification of Algorithm \ref{alg:invexProximal} as described in Section \ref{sub:PnP} \footnote{We used Noise2Void implementation at \url{https://github.com/juglab/n2v}} \cite{krull2019noise2void}. The same deconvolution problem as in Experiment 1 is used. The interesting aspect of this scenario is that Algorithm \ref{alg:invexProximal} has a proximal step in Line 5 that allows to compare between regularizers (invex and convex) while using neural networks in Line 6 (see Algorithm 3 in Appendix \ref{app:PnP} of Supplemental material). Noise2Void is trained by randomly extracting patches of size $64\times 64$ pixels from the training images where zero-mean white Gaussian noise was added for $SNR = 20,30$dB. Data augmentation on the training dataset is used, by rotating each image three times by $90°$ and also added all mirrored versions. The learning rate is fixed as $0.0004$. 
	
	\begin{table}[ht]
		\renewcommand\arraystretch{1.2}
		\centering
		\caption{Performance comparison, in terms of PSNR (dB), where the best and least efficient among invex functions is highlighted in boldface and underscore, respectively.\vspace{1em}}
		\begin{subtable}{1\textwidth}
			\resizebox{1\textwidth}{!}{ \renewcommand{\arraystretch}{1.3}
				\begin{tabular}{P{0.5cm} P{1cm} P{1cm} P{1cm} P{1cm} P{1cm} P{1.2cm} |P{2.0cm}| P{2.5cm} |P{1.8cm} P{1.8cm} P{1.8cm}}
					\hline
					\multicolumn{7}{P{7.5cm}|}{(Experiment 1) Algorithm \ref{alg:invexProximal}, $p=0.5$ for Eq. \eqref{fun1}.} & FISTA \cite{beck2009fast} & ReconNet \cite{kulkarni2016reconnet} & TVAL3 \cite{li2013efficient} & SCAD \cite{fan2001variable} & MCP \cite{zhang2010nearly} \\
					\hline
					\multicolumn{2}{P{1.3cm}}{SNR} & Eq. \eqref{fun1} & Eq. \eqref{fun2} & Eq. \eqref{fun3} & Eq. \eqref{fun4} & Eq. \eqref{fun5} & $\ell_{1}$-norm & & & & \\
					\hline
					\multicolumn{2}{P{1.3cm}}{\centering $\infty$} & \textbf{33.40} & 31.25 & 31.93 & $\underline{30.00}$ & 32.65 & 29.97 & 27.01 & 28.77 & 30.55 & 31.30 \\
					\multicolumn{2}{P{1.3cm}}{\centering $20$dB} & \textbf{24.60} & 22.83 & 23.39 & $\underline{22.00}$ & 23.98 & 21.80 & 19.99 & 20.49 & 22.60 & 23.01 \\ 
					\multicolumn{2}{P{1.3cm}}{\centering $30$dB} & \textbf{27.61} & 26.56 & 26.90 & $\underline{26.00}$ & 27.25 & 24.91 & 22.01 & 23.99 & 26.10 & 26.77 \\
					\hline
				\end{tabular}
			}
		\end{subtable}
		\medskip
		\medskip
		
		\begin{subtable}{0.48\textwidth}
			\centering
			\resizebox{1\textwidth}{!}{ \renewcommand{\arraystretch}{1.3}
				\begin{tabular}{P{0.5cm} P{1cm} P{1cm} P{1cm} P{1cm} P{1.2cm} P{1.2cm}}
					\hline
					\multicolumn{7}{P{7.5cm}}{(Experiment 2) Algorithm \ref{alg:invexPnP}, $p=0.8$ for Eq. \eqref{fun1}.} \\
					\hline
					SNR & Eq. \eqref{fun1} & Eq. \eqref{fun2} & Eq. \eqref{fun3} & Eq. \eqref{fun4} & Eq. \eqref{fun5} & $\ell_{1}$-norm \\
					\hline
					\centering $\infty$ & \textbf{34.51} & 32.37 & 33.06 & $\underline{31.40}$ & 33.76 & 31.10 \\
					\centering $20$dB & \textbf{25.55} & 23.92 & 24.44 & $\underline{23.00}$ & 24.98 & 22.95 \\
					\centering $30$dB & \textbf{28.30} & 26.87 & 27.33 & $\underline{26.05}$ & 27.80 & 26.00 \\
					\hline
				\end{tabular}
			}
		\end{subtable}
		\hfil
		\begin{subtable}{0.48\textwidth}
			\centering
			\resizebox{1\textwidth}{!}{ \renewcommand{\arraystretch}{1.45}
				\begin{tabular}{P{1.5cm} P{1cm} P{1cm} P{1.2cm} P{1.5cm} P{1.5cm}}
					\hline
					\multicolumn{4}{P{6cm}|}{(Denoising experiment) Algorithm \ref{alg:denoising}, $p=0.5$ for Eq. \eqref{fun1}} & BM3D \cite{dabov2007image} & Noise2Void \cite{krull2019noise2void} \\
					\hline
					Metric & Eq. \eqref{fun1} & Eq. \eqref{fun3} & Eq. \eqref{fun5} & $\ell_{1}$-norm & \\
					\hline
					SNR (dB) & \textbf{49.40} & $\underline{43.85}$ & 46.46 & 41.52 & 39.43 \\
					SSIM & \textbf{0.886} & $\underline{0.872}$ & 0.876 & 0.869 & 0.853 \\
					\hline
				\end{tabular}
			}
		\end{subtable}
		\medskip
		\medskip
		
		\begin{subtable}{0.48\textwidth}
			\centering
			\resizebox{1\textwidth}{!}{ \renewcommand{\arraystretch}{1.3}
				\begin{tabular}{P{0.5cm} P{1cm} P{1cm} P{1cm} P{1cm} P{1.2cm} P{1.2cm} | P{1.2cm} | P{1.2cm}}
					\hline
					\multicolumn{9}{P{10.5cm}}{(Experiment 3) Algorithm \ref{alg:unrolling} - unfolded LISTA. $p=0.85$ for Eq. \eqref{fun1}} \\
					\hline
					SNR & $m/n$& Eq. \eqref{fun1} & Eq. \eqref{fun2} & Eq. \eqref{fun3} & Eq. \eqref{fun4} & Eq. \eqref{fun5} & $\ell_{1}$-norm \cite{chen2018theoretical} & ReconNet \cite{kulkarni2016reconnet} \\
					\hline
					\centering \multirow{3}{*}{$\infty$} & \centering 0.2 \vspace{0.1em} \hrule& \textbf{31.32} & 29.20 & 29.87 & $\underline{28.56}$ & 30.58 & 27.95 & 26.59 \\ 
					& \centering 0.4 \vspace{0.1em} \hrule& \textbf{36.10} & 33.50 & 34.34 & $\underline{32.75}$ & 35.20 & 32.01 & 31.86 \\
					& \centering 0.6 & \textbf{41.27} & 37.81 & 38.90& $\underline{36.09}$ & 40.05 & 35.82 & 34.42 \\ 
					\hline
					\centering \multirow{3}{*}{20dB} & \centering 0.2 \vspace{0.1em} \hrule& \textbf{26.00} & 24.45 & 24.94 & $\underline{23.97}$ & 25.01 & 23.52 & 22.00 \\ 
					& \centering 0.4 \vspace{0.1em} \hrule& \textbf{32.67} & 30.64 & 31.32 & $\underline{30.02}$ & 32.29 & 29.43 & 28.24 \\ 
					& \centering 0.6 & \textbf{34.38} & 33.00 & 33.28 & $\underline{32.94}$ & 33.64 & 32.60 & 30.20 \\ 
					\hline
					\centering \multirow{3}{*}{30dB} &\centering 0.2 \vspace{0.1em} \hrule& \textbf{27.65} & 26.20 & 26.66 & $\underline{25.75}$ & 27.15 & 25.32 & 23.64 \\ 
					& \centering 0.4 \vspace{0.1em} \hrule& \textbf{34.33} & 31.89 & 32.66 & $\underline{31.02}$ & 33.47 & 30.46 & 29.88 \\ 
					& \centering 0.6 & \textbf{37.03} & 34.84 & 35.54 & $\underline{34.17}$ & 36.27 & 33.53 & 31.71 \\ 
					\hline
				\end{tabular}
			}
		\end{subtable}
		\hfil
		\begin{subtable}{0.48\textwidth}
			\centering
			\resizebox{1\textwidth}{!}{ \renewcommand{\arraystretch}{1.12}
				\begin{tabular}{P{0.5cm} P{1cm} P{1cm} P{1cm} P{1cm} P{1.2cm} | P{1.2cm} | P{1.2cm}}
					\hline
					\multicolumn{8}{P{10.5cm}}{(Experiment 3) Algorithm \ref{alg:unrolling} - unfolded ISTA-Net. $p=0.85$ for Eq. \eqref{fun1}} \\
					\hline
					SNR & $m/n$& Eq. \eqref{fun1} & Eq. \eqref{fun2} & Eq. \eqref{fun3} & Eq. \eqref{fun4} & Eq. \eqref{fun5} & $\ell_{1}$-norm \cite{zhang2018ista} \\
					\hline
					\centering \multirow{3}{*}{$\infty$} & \centering 0.2 \vspace{0.1em} \hrule& \textbf{32.50} & 30.15 & 30.89 & $\underline{29.04}$ & 31.67 & 28.77 \\
					& \centering 0.4 \vspace{0.1em} \hrule& \textbf{38.33} & 35.72 & 36.55 & $\underline{34.92}$ & 37.41 & 34.17 \\
					& \centering 0.6 & \textbf{43.61} & 40.07 & 41.18& $\underline{39.02}$ & 42.36 & 38.02 \\ 
					\hline
					\centering \multirow{3}{*}{20dB} & \centering 0.2 \vspace{0.1em} \hrule& \textbf{28.29} & 26.22 & 26.87 & $\underline{25.60}$ & 27.56 & 25.01 \\ 
					& \centering 0.4 \vspace{0.1em} \hrule& \textbf{33.96} & 32.11 & 32.71 & $\underline{31.55}$ & 33.32 & 31.00 \\ 
					& \centering 0.6 & \textbf{35.77} & 34.68 & 35.03 & $\underline{34.33}$ & 35.39 & 33.99 \\ 
					\hline
					\centering \multirow{3}{*}{30dB} &\centering 0.2 \vspace{0.1em} \hrule& \textbf{29.34} & 28.30 & 28.63 & $\underline{27.97}$ & 28.98 & 27.65 \\ 
					& \centering 0.4 \vspace{0.1em} \hrule& \textbf{35.41} & 33.33 & 33.99 & $\underline{32.69}$ & 34.68 & 32.08 \\ 
					& \centering 0.6 & \textbf{38.95} & 36.25 & 37.10 & $\underline{35.43}$ & 38.00 & 34.65 \\ 
					\hline
				\end{tabular}
			}
		\end{subtable}
		\label{tab:globalResults}
	\end{table}

	\textbf{Experiment 3} compares the invex regularizers but under the unrolling framework as described in Section \ref{unrolling}. The gold standard convex regularizations to compare with are the learned iterative shrinkage and thresholding algorithm (LISTA) \cite{liu2019alista}, and the Interpretable optimization-inspired deep network (ISTA-Net)\cite{zhang2018ista}. Also, to comparing with convolutional neural networks methodologies, the non-iterative reconstruction methodology ReconNet \cite{kulkarni2016reconnet} is used. We follow the existing setting for LISTA in \cite{chen2018theoretical}\footnote{We used the implementation from \cite{chen2018theoretical} at \url{https://github.com/VITA-Group/LISTA-CPSS}}, and for ISTA-Net in \cite{liu2019alista}. For the training stage we extract $10000$ patches $\boldsymbol{b}\in \mathbb{R}^{16\times 16}$ at random positions of each image, with all means removed. We then learn a dictionary $\boldsymbol{D}\in \mathbb{R}^{256\times 512}$ from the extracted patches, using the same strategy as in \cite{chen2018theoretical}. Gaussian i.i.d sensing matrices $\boldsymbol{\Phi}\in \mathbb{R}^{m\times 256}$ are created from the standard Gaussian distribution, $\boldsymbol{\Phi}[i,j]\sim \mathcal{N}(0,1/m)$ and then normalize its columns to have the unit $\ell_{2}$-norm, where $m$ is selected such that $\frac{m}{256}=0.2,0.4,0.6$. The matrix $\boldsymbol{H}$ is built as $\boldsymbol{H}=\boldsymbol{\Phi}\boldsymbol{\Psi}$ with $T=16$ (number of layers). We follow the same two-step strategy in \cite{chen2018theoretical} to train a recurrent neural network. First, perform a layer-wise pre-training solving Eq. \eqref{eq:problem4} for each extracted patch $\boldsymbol{b}$ by fixing $\boldsymbol{H}=\boldsymbol{\Psi}$. Second, append a learnable fully-connected layer at the end of the network structure, initialized by $\boldsymbol{\Psi}$. Then, perform an end-to-end training solving Eq. \eqref{eq:problem4} where $\boldsymbol{H}$ in this case is learnt by updating the initial matrix $\boldsymbol{\Psi}$. For each testing image, we divide it into non-overlapping $16\times 16$ patches. When $g(\boldsymbol{x})$ is the the $\ell_{1}$-norm, we recover~\cite{chen2018theoretical}. 
	
	In the case of ISTA-Net, and ReconNet, for their learning stage we extract patches from the training image set of size $33\times 33$. Gaussian i.i.d sensing matrices $\boldsymbol{\Phi}\in \mathbb{R}^{m\times 1089}$ are created with $\ell_{2}$-normalized columns as for LISTA, where $m$ is selected such that $\frac{m}{1089}=0.2,0.4,0.6$. The optimizer employed was Adam algorithm and a learning rate $1\times 10^{-4}$ for 200 and 512 epochs for ISTA-Net and ReconNet respectively, with a batch size of $64$ for both networks. For ISTA-Net $T=16$ (number of unrolled iterations). We recall that when $g(\boldsymbol{x})$ is the the $\ell_{1}$-norm in ISTA-Net, we recover~\cite{zhang2018ista}.\vspace{-0.5em}

	\subsection{ Image Denoising Experiment}
	Two image datasets, which we merge ($80$ images in total), are used for this experiment comes from a neutron image formation phenomenon\footnote{Acquired with the ISIS Neutron and Moun Source system at Harwell Science and Innovation Campus.}. These type of images contain the neutron attenuation properties of the object which helps analyze material structure. Performance is assessed by averaging along all the images the experimental SNR in dB given by $SNR = 20\log\left (\frac{\lVert \boldsymbol{z} \rVert_{2}}{\lVert \hat{\boldsymbol{z}}-\boldsymbol{z} \rVert_{2}}\right)$, where $\boldsymbol{z}$ and $\hat{\boldsymbol{z}}$ stand for the noisy and the denoised image, respectively, and the structural similarity index measure (SSIM) computed between $\boldsymbol{z}$ and $\hat{\boldsymbol{z}}$. Taking advantage of results observed from previous experiments, we compare the top three regularizers in Eqs. \eqref{fun1}, \eqref{fun3}, and \eqref{fun5} with two state-of-the-art denoising techniques including the block-matching and 3-D filtering (BM3D) \cite{dabov2007image} using $\ell_{1}$-norm regularizer and the deep learning technique Noise2Void (trained as in Experiment 2) \cite{krull2019noise2void}. We follow the two-step denoising procedure described in Section \ref{sub:denoising}. In the first step, the transform domain is built using PCA as in \cite{cai2014data}. To build this transform we extract patches of $16\times 16$ from the noisy image that are then used to adaptively construct a tight frame (nearly orthogonal matrix) tailored to the given noisy data \footnote{We used implementation at \url{https://www.math.hkust.edu.hk/~jfcai/}.}. Results are summarized in Table \ref{tab:globalResults}. We report examples of denoised images obtained by Eqs. \eqref{fun1}, \eqref{fun3}, \eqref{fun5}, BM3D, and Noise2Void are illustrated in Appendix \ref{app:denoising} of supplementary material, along with the algorithm used for the invex regularizers to denoise these images.\vspace{-1.em}
	
	%

	\section{Discussion, Limitations and Conclusion}\vspace{-0.8em}
	Application advancement of invex theory has paused for decades due to the lack of practical examples, which has caused a significantly reduced interest in invexity research. To address this issue, we present for the first time a list of invex regularizers for image reconstruction applications, and formulate corresponding optimization problems. Particularly, for image compressive sensing, we improve three advanced imaging techniques using the listed functions in Table \ref{tab:list} as invex regularizers. We present their solution algorithms and develop theoretical guarantees on their convergence to global minimum. We also conducted various image compressive sensing and denoising experiments to demonstrate the effectiveness of invex regularizers under practical scenarios that are non-ideal with noisy data observed and RIP condition not guaranteed. Significant benefit of using invex regularizers have been proved from both theoretical and empirical aspects. In fact, Table \ref{tab:globalResults} and theoretical results in Section \ref{sec:invexImag} revive the potential of exploring invex theory in practical applications.
	
	The numerical results presented in Table \ref{tab:globalResults} confirm performance improvement by using invex regularizers over the $\ell_{1}$-norm-based methods (e.g FISTA, TVAL3) in unexplored scenarios. These tables and theoretical results in Section \ref{sec:invexImag} revive the potential of exploring invex theory in practical applications. The best result is obtained with Eq. \eqref{fun1}, and Eq. \eqref{fun4} is the least efficient. The intuition behind the superiority of Eq. \eqref{fun1} comes from the possibility of adjusting the value of $p$ in data-dependent manner \cite{wu2013improved}. This means that when the images are strictly sparse, and the noise is relatively low, a small value of $p$ should be used. Conversely, when images are non-strictly sparse and/or the noise is relatively high, a larger value of $p$ tend to yield better performance (which seems to be the case for the selected image datasets). We believe that the remaining invex, SCAD, and MCP regularizers have a lesser performance than Eq. \eqref{fun1} as they do not have the flexibility of adjustment to the sparsity of the data. In fact, Eq. \eqref{fun4} shows the poorest performance because in the proof of Theorem \ref{theo:ourCS}, we theoretically guarantee that Eq. \eqref{fun4} cannot sparsify all images. Therefore, this analysis leads to the conclusion that the invex function Eq. \eqref{fun1} offers the best performance for the metrics concerns and the imaging problems studied here.
		
	Although, we have presented theoretical results with global optima using invexity for some of most important and successful image reconstruction techniques, we highlight several limitations of our analysis. Specifically, we focused on reconstructed the image of interest in an ideal scenario, that is, without the present of noise (Theorem \ref{theo:ourCS}). Additionally, we have limited our numerical results to tasks like denoising, and deconvolution. And, the convergence guarantees for the plug-and-play result only ensures a close estimate of the solution (Lemma \ref{theo:PnP}). Therefore, we see there are a number of future directions this research can be taken further improving the results even further. One aspect is to explore avenues for improving convergence guarantees to global optima the plug-and-play framework. Another direction is the study of inclusion of noise in the analysis of imaging applications, which may be an enabler to improve downstream tasks like invex robust image reconstruction. Finally, we feel that the application domains for invex functions can go well beyond denoising, and deconvolution imaging problems, especially around deep learning research, which can improved a number of downstream applications.\vspace{-1em}
	
	\section*{Broader Impact}\vspace{-1em}
	We believe that the presented mathematical and empirical analysis over the studied regularizers has the potential to unlock the benefits of invexity for further applications in signal and image processing. This may be an enabler to improve downstream tasks like deep learning for imaging, and to provide more robust image reconstruction algorithms.\vspace{-1em}
	
	\section*{Acknowledgments}\vspace{-1em}
	This work was partially supported by the Facilities Funding from Science and Technology Facilities Council (STFC) of UKRI, and Wave 1 of the UKRI Strategic Priorities Fund under the EPSRC grant EP/T001569/1, particularly the ``AI for Science" theme within that grant, by the Alan Turing Institute.
	
	{
		\small
		\bibliographystyle{unsrt}
		\bibliography{sample}
	}

	\newpage
	\appendix

	\section{Proof of Lemma \ref{theo:invexProof}}
	\label{app:invexProof}
	In this proof we seek to guarantee that the list of functions in Table \ref{tab:list} are invex. We point out that, since the regularizers in Table \ref{tab:list} is the sum of a scalar function applied to each entry of a vector, then it is enough to analyze the scalar function to determine the invexity of the regularizer.
	
	\paragraph{Eq. \eqref{fun1}.}
	\begin{proof}
		Take $r_{\epsilon}(w) = \left(\lvert w \rvert + \epsilon \right)^{p}, \forall w\in \mathbb{R}$, for $p\in (0,1)$ and $\epsilon\geq \left(p(1-p)\right)^{\frac{1}{2-p}}$. The need to add the constant $\epsilon$ it is to formally satisfy the Lipschitz continuous condition required to be invex according to Definition \ref{def:invex}. Observe that if $w>0$ then we have that $\partial r_{\epsilon}(w)=\left \lbrace \frac{p}{\left(\lvert w \rvert + \epsilon \right)^{1-p}}\right \rbrace$, which means that $0\not \in \partial r_{\epsilon}(w)$. Conversely, if $w<0$ then $\partial r_{\epsilon}(w)=\left \lbrace \frac{-p}{\left(\lvert w \rvert + \epsilon \right)^{1-p}}\right \rbrace$, leading to $0\not \in \partial r_{\epsilon}(w)$. Lets examinate $w^{*}=0$. Note that 
		\begin{align}
		\lim_{w\rightarrow 0^{+}} r_{\epsilon}^{\prime}(w) = \lim_{w\rightarrow 0^{+}} \frac{p}{\left(\lvert w \rvert + \epsilon \right)^{1-p}} = \frac{p}{\epsilon^{1-p}},
		\label{eq:rightLimit}
		\end{align}
		and that 
		\begin{align}
		\lim_{w\rightarrow 0^{-}} r_{\epsilon}^{\prime}(w) = \lim_{w\rightarrow 0^{-}} \frac{-p}{\left(\lvert w \rvert + \epsilon \right)^{1-p}} = \frac{-p}{\epsilon^{1-p}}.
		\label{eq:leftLimit}
		\end{align}
		Additionally, since $r_{\epsilon}(w)$ is a Lipschitz continuous function, then appealing to Theorem \ref{theo:auxDerivative} we have that $\partial r_{\epsilon}(w^{*}=0) = \text{ conv }\left \lbrace \frac{-p}{\epsilon^{1-p}}, \frac{p}{\epsilon^{1-p}}\right \rbrace = \left \lbrack \frac{-p}{\epsilon^{1-p}}, \frac{p}{\epsilon^{1-p}} \right\rbrack$. This means that $0\in \partial r_{\epsilon}(0)$. Further, given the fact that $r_{\epsilon}(0)\leq r_{\epsilon}(w)$ for all $w\in \mathbb{R}$, then $w^{*} = 0$ is a global minimizer of $r_{\epsilon}$. Therefore, the function $r_{\epsilon}$ is invex.
	\end{proof}
	
	\paragraph{Eq. \eqref{fun2}}
	\begin{proof}
		Take $r(w) = \log(1+\lvert w\rvert)$. Observe that if $w>0$ then we have that $\partial r(w)=\left \lbrace \frac{1}{1 + \lvert w \rvert}\right \rbrace$, which means that $0\not \in \partial r(w)$. Conversely, if $w<0$ then $\partial r(w)=\left \lbrace \frac{-1}{1 + \lvert w \rvert}\right \rbrace$, leading to $0\not \in \partial r(w)$. Lets examinate $w^{*}=0$. Note that 
		\begin{align}
		\lim_{w\rightarrow 0^{+}} r^{\prime}(w) = \lim_{w\rightarrow 0^{+}} \frac{1}{1 + \lvert w \rvert} = 1,
		\label{eq:rightLimit1}
		\end{align}
		and that 
		\begin{align}
		\lim_{w\rightarrow 0^{-}} r^{\prime}(w) = \lim_{w\rightarrow 0^{-}} \frac{-1}{1+ \lvert w \rvert} = -1.
		\label{eq:leftLimit1}
		\end{align}
		Additionally, since $r(w)$ is a Lipschitz continuous function, then appealing to Theorem \ref{theo:auxDerivative} we have that $\partial r(w^{*}=0) = \text{ conv }\left \lbrace -1,1\right \rbrace = \left \lbrack -1,1 \right\rbrack$. This means that $0\in \partial r(0)$. Further, given the fact that $r(0)\leq r(w)$ for all $w\in \mathbb{R}$, then $w^{*} = 0$ is a global minimizer of $r(w)$. Therefore, the function $r(w)$ is invex.
	\end{proof}
	
	\paragraph{Eq. \eqref{fun3}}
	\begin{proof}
		Take $r(w) = \frac{\lvert w \rvert}{2+2\lvert w \rvert}$. Observe that if $w>0$ then we have that $\partial r(w)=\left \lbrace \frac{1}{2(1 + \lvert w \rvert)^{2}}\right \rbrace$, which means that $0\not \in \partial r(w)$. Conversely, if $w<0$ then $\partial r(w)=\left \lbrace \frac{-1}{2(1 + \lvert w \rvert)^{2}}\right \rbrace$, leading to $0\not \in \partial r(w)$. Lets examinate $w^{*}=0$. Note that 
		\begin{align}
		\lim_{w\rightarrow 0^{+}} r^{\prime}(w) = \lim_{w\rightarrow 0^{+}} \frac{1}{2(1 + \lvert w \rvert)^{2}} = \frac{1}{2},
		\label{eq:rightLimit2}
		\end{align}
		and that 
		\begin{align}
		\lim_{w\rightarrow 0^{-}} r^{\prime}(w) = \lim_{w\rightarrow 0^{-}} \frac{-1}{2(1 + \lvert w \rvert)^{2}} = -\frac{1}{2}.
		\label{eq:leftLimit2}
		\end{align}
		Additionally, since $r(w)$ is a Lipschitz continuous function, then appealing to Theorem \ref{theo:auxDerivative} we have that $\partial r(w^{*}=0) = \text{ conv }\left \lbrace -\frac{1}{2},\frac{1}{2}\right \rbrace = \left \lbrack -\frac{1}{2},\frac{1}{2} \right\rbrack$. This means that $0\in \partial r(0)$. Further, given the fact that $r(0)\leq r(w)$ for all $w\in \mathbb{R}$, then $w^{*} = 0$ is a global minimizer of $r(w)$. Therefore, the function $r(w)$ is invex.
	\end{proof}
	
	\paragraph{Eq. \eqref{fun4}}
	\begin{proof}
		Consider $r(w) = \frac{w^{2}}{1+w^{2}}$. Observe that $\partial r(w) = \left\lbrace \frac{2w}{(1+w^{2})^{2}}\right \rbrace$, which means $r(w)$ is continuously differentiable. Then, it is clear that $w=0$ is the only point that satisfies $0\in \partial r(0)$. In addition, the value $r(w=0)$ is the global minimum of $r(w)$. Thus, since the only stationary point of $r(w)$ is a global minimizer, then $r(w)$ is invex.
	\end{proof}
	
	\paragraph{Eq. \eqref{fun5}}
	\begin{proof}
		Take $r(w) = \log(1+\lvert w \rvert) - \frac{\lvert w \rvert}{2 + 2\lvert w \rvert}$. Observe that if $w>0$ then we have that $\partial r(w)=\left \lbrace \frac{1}{2(1+\lvert w \rvert)^{2}} + \frac{w}{(1+\lvert w \rvert)^{2}} \right \rbrace$, which means that $0\not \in \partial r(w)$. Conversely, if $w<0$ then $\partial r(w)=\left \lbrace \frac{-1}{2(1+\lvert w \rvert)^{2}} + \frac{w}{(1+\lvert w \rvert)^{2}}\right \rbrace$, leading to $0\not \in \partial r(w)$. Lets examinate $w^{*}=0$. Note that 
		\begin{align}
		\lim_{w\rightarrow 0^{+}} r^{\prime}(w) = \lim_{w\rightarrow 0^{+}} \frac{1}{2(1+\lvert w \rvert)^{2}} + \frac{w}{(1+\lvert w \rvert)^{2}} = \frac{1}{2},
		\label{eq:rightLimit3}
		\end{align}
		and that 
		\begin{align}
		\lim_{w\rightarrow 0^{-}} r^{\prime}(w) = \lim_{w\rightarrow 0^{-}} \frac{-1}{2(1+\lvert w \rvert)^{2}} + \frac{w}{(1+\lvert w \rvert)^{2}} = -\frac{1}{2}.
		\label{eq:leftLimit4}
		\end{align}
		Additionally, since $r(w)$ is a Lipschitz continuous function, then appealing to Theorem \ref{theo:auxDerivative} we have that $\partial r(w^{*}=0) = \text{ conv }\left \lbrace -\frac{1}{2},\frac{1}{2}\right \rbrace = \left \lbrack -\frac{1}{2},\frac{1}{2} \right\rbrack$. This means that $0\in \partial r(0)$. Further, given the fact that $r(0)\leq r(w)$ for all $w\in \mathbb{R}$, then $w^{*} = 0$ is a global minimizer of $r(w)$. Therefore, the function $r(w)$ is invex.
	\end{proof}
	
	\subsection{Additional Discussion on Invex Regularizers}
	\label{app:discussionRegu}
	To address sub-optimal limitations of convex regularizers, non-convex mappings have been proposed. For instance, the Smoothly Clipped Absolute Deviation (SCAD) \cite{fan2001variable}, and Minimax Concave Penalty (MCP) \cite{zhang2010nearly}. However, a recent survey in imaging \cite{wen2018survey}, which compared the performance of several regularizers including SCAD and MCP for a number of imaging, concludes that Eq. \eqref{fun1} shows higher performance than SCAD and MCP because the value of $p$ can be adjusted in data-dependent manner. This means that when the images are strictly sparse, and the noise is relatively low, a small value of $p$ should be used. Conversely, when images are non-strictly sparse and/or the noise is relatively high, a larger value of $p$ tend to yield better performance. Furthermore, in the context of invexity, we highlight that SCAD and MCP are non-invex regularizer because they reach a maximum value, which makes the first derivative zero in non-minimizer values leading to its non-invexity (see Theorem \ref{theo:invexProof}).
	
	On the other hand, in the case of minimax-concave-type of regularizers, we present a new function  in our manuscript (Eq. \eqref{fun5}). From Eq. \eqref{fun5} it is clear we are subtracting $g_{1}(\boldsymbol{x})=\sum_{i=1}^{n}\log(1+\lvert \boldsymbol{x}[i] \rvert)$, and $g_{2}(\boldsymbol{x})=\sum_{i=1}^{n}\frac{\lvert \boldsymbol{x}[i] \rvert}{2 + 2\lvert \boldsymbol{x}[i] \rvert}$ (selected due to results in \cite{wu2019improved}). We propose to study regularizer in Eq. \eqref{fun5}, that is $g_{1}-g_{2}$, for three reasons. First, because $g_{1}$, $g_{2}$, and $g_{1}-g_{2}$ are invex, as stated in Lemma 1, and all of them can achieve global optima for the scenarios studied in the paper. Second, to the best of our knowledge, there is no evidence that subtracting two convex penalties (current proposal in the minimax-concave literature) produces another convex regularizer (if exists). Therefore, we present Eq. \eqref{fun5} to show that at least this is possible in the invex case, as stated in Section \ref{sec:InvexRegul}. 
	
	Finally, we point out that the performance of invex regularizers in Eqs. \eqref{fun2}, \eqref{fun3}, \eqref{fun4}, and \eqref{fun5} can be justified under the framework of re-weighted $\ell_{1}$-norm minimization (see \cite{candes2008enhancing}), which enhances the performance of just $\ell_{1}$-norm minimization.
	
	\section{Proof of Lemma \ref{theo:invexComposited}}
	\label{app:invexComposited}
	\begin{proof}
		To prove this theorem, we show that for each $\boldsymbol{x}\in \mathbb{R}^{n}$ such that $\boldsymbol{0} \in \partial h(\boldsymbol{x})$ where $h(\boldsymbol{x}) = f\left(\boldsymbol{H}\boldsymbol{x} -\boldsymbol{v}\right)$ is a global minimizer. Observe that 
		\begin{align}
		\partial h(\boldsymbol{x}) = \left\{ \nabla h(\boldsymbol{x}) \right\} = \left\{\boldsymbol{H}^{T}\nabla f(\boldsymbol{H}\boldsymbol{x}-\boldsymbol{v})\right\}.
		\label{eq:gradient1}
		\end{align}
		Take $\boldsymbol{x}^{*}\in \mathbb{R}^{n}$ such that $\boldsymbol{0}\in \partial h(\boldsymbol{x}^{*})$, then since $\boldsymbol{H}$ is a full row-rank matrix (equivalently $\boldsymbol{H}^{T}$ full col-rank matrix) from Eq. \eqref{eq:gradient1} we have
		\begin{align}
		\nabla h(\boldsymbol{x}^{*}) = \boldsymbol{H}^{T}\nabla f(\boldsymbol{H}\boldsymbol{x}^{*}-\boldsymbol{v}) = \boldsymbol{0}\leftrightarrow \nabla f(\boldsymbol{H}\boldsymbol{x}^{*}-\boldsymbol{v}) =\boldsymbol{0}.
		\label{eq:problem2}
		\end{align}
		The above equation means that each stationary point of $h(\boldsymbol{x})$ is found through the stationary points of $f(\boldsymbol{x})$. Thus, since $f$ is invex then $\boldsymbol{H}\boldsymbol{x}^{*}-\boldsymbol{v}$ is a global minimizer of $f$ i.e. $h$ is invex.
	\end{proof}
	
	\section{Proof of Theorem \ref{theo:proximalProof}}
	\label{app:proximalProof}
	In this appendix we seek to guarantee that the proximal operator of the functions in Table \ref{tab:list} are invex. We point out that, since the proximal of the regularizers in Table \ref{tab:list} is the sum of a scalar function applied to each entry of a vector, then it is enough to analyze the scalar function to determine the invexity of the proximal.
	
	\subsection{Invexity proofs of the proximal operators}
	In the following we provide the proof for the first statement in Theorem \ref{theo:proximalProof}.
	
	\paragraph{Eq. \eqref{fun1}}
	\begin{proof}
		Let $h(w)$ be a function defined, for $p\in (0,1)$, as
		\begin{align}
		h(w) = (\lvert w \rvert + \epsilon)^{p} + \frac{1}{2}(w-u)^{2},
		\end{align}
		for fixed $u\in \mathbb{R}$, and $\epsilon\geq \left(p(1-p)\right)^{\frac{1}{2-p}}$. Then, we seek to show that the second derivate of $h(w)$ with respect to $w$ for $w\not = 0$ is non-negative. Observe that,
		\begin{align}
		h^{\prime \prime}(w) = \frac{p(p-1)}{(\lvert w \rvert + \epsilon )^{2-p}} + 1.
		\label{eq:deri1}
		\end{align}
		From Eq. \eqref{eq:deri1} we have that $(\lvert w \rvert + \epsilon )^{2-p}$ is a positive increasing function since $2-p>1$. This implies that to show $h^{\prime \prime}(w)$ is non-negative for all $w\in \mathbb{R}$ we need to analyze only when $w=0$. Therefore, $\frac{p(p-1)}{(\lvert w \rvert + \epsilon )^{2-p}} \in [-1,0)$ for all $w\in \mathbb{R}$, because $p(p-1)<0$ and $\epsilon \geq \left(p(1-p)\right)^{\frac{1}{2-p}}$. Thus, $h^{\prime \prime}(w)$ is non-negative, leading to the invexity of $h(w)$ (i.e. $h^{\prime \prime}(w)$ positive implies convexity).
	\end{proof}
	
	\paragraph{Eq. \eqref{fun2}}
	\begin{proof}
		Take $h(w) = \log(1+\lvert w\rvert) + \frac{1}{2}(w-u)^{2}$ for fixed $u\in \mathbb{R}$. Observe that the second derivative of $h(w)$, for $w\not =0$ is given by
		\begin{align}
		h^{\prime \prime}(w) = \frac{-1}{(1 + \lvert w \rvert)^{2}} + 1.
		\end{align}
		Then, since $(1 + \lvert w \rvert)^{2}\geq 1$ for all $w$, this implies that $\frac{-1}{(1 + \lvert w \rvert)^{2}} \in [-1,0)$. Thus, $h^{\prime \prime}(w)$ is non-negative, leading to the invexity of $h(w)$.
	\end{proof}
	
	\paragraph{Eq. \eqref{fun3}}
	\begin{proof}
		Take $h(w) = \frac{\lvert w \rvert}{2 + 2\lvert w \rvert} + \frac{1}{2}(w-u)^{2}$ for fixed $u\in \mathbb{R}$. We will use the same argument as in previous cases. Then, for $w\not =0$ notice that the second derivative of $h(w)$ is given by
		\begin{align}
		h^{\prime \prime}(w) = \frac{-1}{(1+\lvert w \rvert)^{3}} + 1.
		\end{align}
		Then, from the above equation it is clear that $\frac{-1}{(1+\lvert w \rvert)^{3}}\in [-1,0)$ for all $w\in \mathbb{R}$. Thus, $h^{\prime \prime}(w)$ is non-negative, leading to the invexity of $h(w)$.
	\end{proof}
	
	\paragraph{Eq. \eqref{fun4}}
	\begin{proof}
		Take $h(w) = \frac{w^{2}}{1+w^{2}} + \frac{1}{2}(w-u)^{2}$, for fixed $u\in \mathbb{R}$. Then, notice that the second derivative of $h(w)$ is given by
		\begin{align}
		h^{\prime \prime}(w) = \frac{2-6w^{2}}{(1+ w^{2})^{3}} + 1.
		\end{align}
		Then, we show that $s(w)=\frac{2-6w^{2}}{(1+ w^{2})^{3}} \geq -1$, by determining its extreme values. Observe that
		\begin{align}
		s^{\prime}(w) = \frac{24w(w^{2}-1)}{(1+w^{2})^{3}} =0,
		\end{align}
		only when $w=0,1,-1$. It is clear that the maximum value of $s(w)$ is attained when $w=0$, i.e. $s(w)=2$. And, its minimum value is achieved when $w=-1$, that is $s(1)=s(-1)=\frac{-1}{2}$. Thus, since $s(w)\geq -1$ then $h(w)$ is invex. 
	\end{proof}
	
	\paragraph{Eq. \eqref{fun5}}
	\begin{proof}
		Take $h(w) = \log(1+\lvert w \rvert) - \frac{\lvert w \rvert}{2 + 2\lvert w \rvert} + \frac{1}{2}(w-u)^{2}$, for fixed $u\in \mathbb{R}$. Then, for $w\not =0$ notice that the second derivative of $h(w)$ is given by
		\begin{align}
		h^{\prime \prime}(w) = \frac{-\lvert w \rvert}{(1+\lvert w \rvert)^{3}} + 1.
		\end{align}
		Then, from the above equation it is clear that $\frac{-\lvert w \rvert}{(1+\lvert w \rvert)^{3}} \in [-1,1]$, which implies that $h^{\prime \prime}(w)$ is non-negative for any $w$. Thus, $h(w)$ is invex. 
	\end{proof}
	
	\subsection{The resolvent of proximal operator only has global optimizers}
	\begin{proof}
		Now we proof the second part of Theorem \ref{theo:proximalProof}. From the previous analysis on each proximal operator, we have that $h(\boldsymbol{x})$ is an convex (therefore invex) function, then Theorem \ref{theo:optimal_v0} states that any global minimizer $\boldsymbol{y}$ of $h$ satisfies that $\mathbf{0}\in \partial h(\boldsymbol{y})$. This condition implies that $\mathbf{0} \in \partial g(\boldsymbol{y}) + (\boldsymbol{y}-\boldsymbol{v})$, from which we obtain that $\boldsymbol{y} \in ( \partial g + \mathbf{I})^{-1}(\boldsymbol{v})$. Thus, we have that $\text{\textbf{prox}}_{g}(\boldsymbol{v})=(\partial g + \mathbf{I})^{-1}(\boldsymbol{v})$ from which the result holds.
	\end{proof}

	\begin{figure}[h!]
		\centering
		\includegraphics[width=1\linewidth]{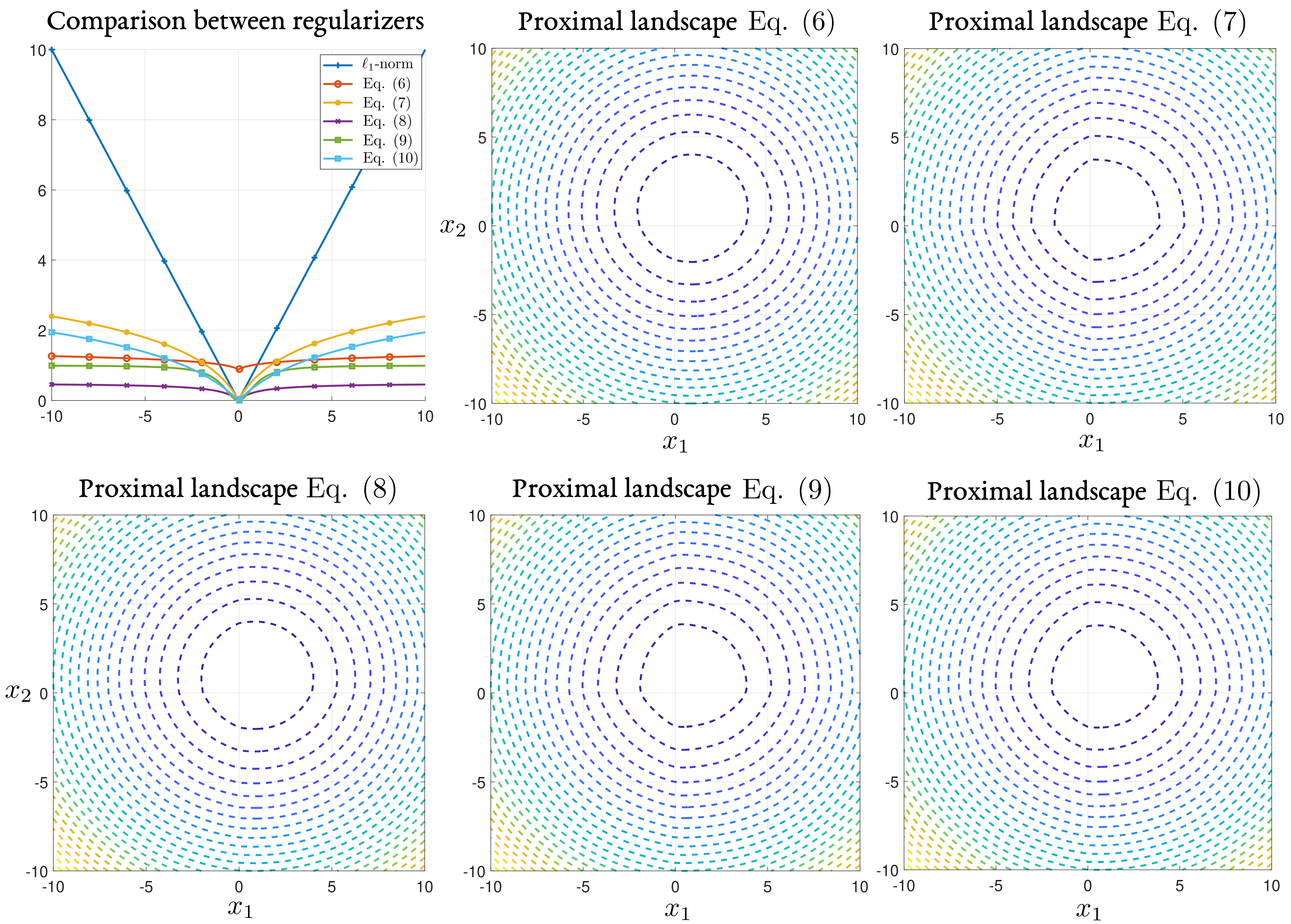}
		\vspace{-0.3em}
		\caption{Here we present a visual comparison between the one-dimensional version of $\ell_{1}$-norm and the invex regularizers in Eqs. \eqref{fun1}, \eqref{fun2}, \eqref{fun3}, \eqref{fun4}, and \eqref{fun5}. For Eq. \eqref{fun1} we select $p=0.5$, and $\epsilon=(p(1-p))^{\frac{1}{2-p}}$. We also report the landscape of function $h(\boldsymbol{x}) = g(\boldsymbol{x}) + \frac{1}{2}\lVert \boldsymbol{x}-\boldsymbol{u} \rVert_{2}^{2}$ where $\boldsymbol{x}$ is a vector of two dimensions $\boldsymbol{x}=[x_{1},x_{2}]^{T}$, $\boldsymbol{u}=[1,1]^{T}$, and for all invex regularizers in Eqs. \eqref{fun1}, \eqref{fun2}, \eqref{fun3}, \eqref{fun4}, and \eqref{fun5}. It is clear that the level curves are concentric convex sets which confirms that $h(\boldsymbol{x})$ is convex (therefore invex), as stated in Theorem \ref{theo:proximalProof}.}
		\label{fig:landscape}
	\end{figure}
	
	\subsection{Numerical Analysis of Proximal}
	In this section we present additional numerical analysis on the proximal of invex regularizers listed in Table \ref{tab:list}. We start by providing a visual comparison between the one-dimensional version of $\ell_{1}$-norm and the invex regularizers in Eqs. \eqref{fun1}, \eqref{fun2}, \eqref{fun3}, \eqref{fun4}, and \eqref{fun5}. This comparison is reported in Fig. \ref{fig:landscape}. From this illustration it is easy to conclude why Eqs. \eqref{fun1}, \eqref{fun2}, \eqref{fun3}, \eqref{fun4}, and \eqref{fun5} are non-convex.

	To complement the comparison between convex and invex regularizers, we present a graphical validation of the theoretical result in Theorem \ref{theo:proximalProof}. To that end, we illustrate also in Fig. \ref{fig:landscape} the landscape of function $h(\boldsymbol{x}) = g(\boldsymbol{x}) + \frac{1}{2}\lVert \boldsymbol{x}-\boldsymbol{u} \rVert_{2}^{2}$ where $\boldsymbol{x}$ is a vector of two dimensions $\boldsymbol{x}=[x_{1},x_{2}]^{T}$, $\boldsymbol{u}=[1,1]^{T}$, with $g(\boldsymbol{x})$ taking the form of all invex regularizers in Eqs. \eqref{fun1}, \eqref{fun2}, \eqref{fun3}, \eqref{fun4}, and \eqref{fun5}. From these results, it is clear that the level curves are concentric convex sets which confirms that $h(\boldsymbol{x})$ is convex (therefore invex), as stated in Theorem \ref{theo:proximalProof}.
	
	Lastly, the running time to compute the proximal of invex regularizers is also an important aspect to compare with its convex competitor i.e. $\ell_{1}$-norm. The reason for this, is because it is desire to improve imaging quality keeping the same computational complexity to obtain it. Therefore, the following Table \ref{tab:timesProx} reports the running time to compute the proximal (in GPU) of Eqs. \eqref{fun1}, \eqref{fun2}, \eqref{fun3}, \eqref{fun4}, and \eqref{fun5} for an image of $2048\times 2048$ pixels. Observe that Table \ref{tab:timesProx} suggests that computing the proximal of the $\ell_{1}$-norm is faster than the proximal of invex regularizers. However, this difference is given in milliseconds making it negligible in practice.
	
	\begin{table}[ht]
		\centering
		\caption{Time to compute the proximal for all invex and convex regularizers, of an image with $2048\times 2048$ pixels. The reported time is the averaged over $256$ trials. For Eq. \eqref{fun1} we select $p=0.5$, and $\epsilon=(p(1-p))^{\frac{1}{2-p}}$.}
		\renewcommand{\arraystretch}{1.2}
		\begin{tabular}{P{0.5cm} P{1cm} P{1cm} P{1cm} P{1cm} P{1cm} P{1.2cm} |P{2.0cm}}
			\hline
			\multicolumn{2}{P{1.3cm}}{} & Eq. \eqref{fun1} & Eq. \eqref{fun2} & Eq. \eqref{fun3} & Eq. \eqref{fun4} & Eq. \eqref{fun5} & $\ell_{1}$-norm \\
			\hline
			\multicolumn{2}{P{1.3cm}}{\centering Time} & $1.47ms$ & $0.63ms$ & $2.8ms$ & $4.7ms$ & $2.4ms$ & $0.66ms$ \\
			\hline
		\end{tabular}
		\label{tab:timesProx}
	\end{table}
	
	\section{Solutions to the Proximal Operator in Eq. \eqref{eq:prox1}}
	\label{app:solProxi}
	
	In this section we present the proximal operator for the functions in Eqs. \eqref{fun1}-\eqref{fun5} summarized in Table \ref{tab:proximals}. In the case of Eq. \eqref{fun1} its proximal operator was calculated in \cite{marjanovic2012l_q}. We recall that the analysis for Eq. \eqref{fun1} is valid with and without the constant $\epsilon$. We prefer to add $\epsilon$ in order to formally satisfy the Lipschitz continuity as in Definition \ref{def:lipschitz}. Moreover, for functions in Eqs. \eqref{fun2}-\eqref{fun5} we present how to estimate their proximal operator them in the following.
	
	\paragraph{Proximal of Eq. \eqref{fun2}}
	Consider $h(w)=\lambda \log(1+\lvert w \rvert) + \frac{1}{2}(w-u)^{2}$ for $\lambda\in (0,1]$, and fixed $u\in \mathbb{R}$. We note first that we only consider $w's$ for which $\text{sign}(w)=\text{sign}(u)$, otherwise $h(w)=\lambda \log(1+\lvert w \rvert) + \frac{1}{2}w^{2}+\lvert u \rvert\lvert w \rvert + \frac{1}{2}u^{2}$ which is clearly minimized at $w=0$. Then, since with $\text{sign}(w)=\text{sign}(u)$ we have $(w-u)^{2} = (\lvert w \rvert-\lvert u \rvert)^{2}$, we replace $u$ with $\lvert u \rvert$ and take $w\geq 0$. As $h(w)$ is differentiable for $w>0$, re-arranging $h^{\prime}(w)=0$ gives
	\begin{align}
	\psi_{\lambda}(w)\delequal\frac{\lambda}{1+w} + w = \lvert u \rvert.
	\label{eq:proxlog}
	\end{align}
	Observe that $\psi^{\prime}_{\lambda}(w)$ is always positive then it means that $\psi_{\lambda}(w)$ is monotonically increasing. Thus, the equation $\psi_{\lambda}(w)=\lvert u \rvert $ has unique solution i.e. at some point the quality holds. Thus, solving $\psi_{\lambda}(w)=\lvert u \rvert$ is equivalent to
	\begin{align}
	w^{2} + (1-\lvert u \rvert) w + \lambda - \lvert u \rvert = 0.
	\label{eq:proxlog1}
	\end{align}
	It is easy to verify that the solution to Eq. \eqref{eq:proxlog1} that returns the minimum value of $h(w)$ is given by $w=\frac{\lvert u \rvert -1 + \sqrt{(\lvert u \rvert + 1)^2 - 4\lambda}}{2}$ when $(\lvert u \rvert+1)^2\geq 4\lambda$, and $0$ otherwise.
	
	\paragraph{Proximal of Eq. \eqref{fun3}}
	Consider $h(w)=\lambda \frac{\lvert w \rvert}{2 + 2\lvert w \rvert} + \frac{1}{2}(w-u)^{2}$ for $\lambda\in (0,1]$, and fixed $u\in \mathbb{R}$. We note first that we only consider $w's$ for which $\text{sign}(w)=\text{sign}(u)$, otherwise $h(w)=\lambda \frac{\lvert w \rvert}{2 + 2\lvert w \rvert} + \frac{1}{2}w^{2}+\lvert u \rvert\lvert w \rvert + \frac{1}{2}u^{2}$ which is clearly minimized at $w=0$. Then, since with $\text{sign}(w)=\text{sign}(u)$ we have $(w-u)^{2} = (\lvert w \rvert-\lvert u \rvert)^{2}$, we replace $u$ with $\lvert u \rvert$ and take $w\geq 0$. As $h(w)$ is differentiable for $w>0$, re-arranging $h^{\prime}(w)=0$ gives
	\begin{align}
	\psi_{\lambda}(w)\delequal\frac{\lambda}{2(1+w)^{2}} + w = \lvert u \rvert.
	\label{eq:proxrat}
	\end{align}
	Observe that $\psi^{\prime}_{\lambda}(w)$ is always positive then it means that $\psi_{\lambda}(w)$ is monotonically increasing. Thus, the equation $\psi_{\lambda}(w)=\lvert u \rvert $ has unique solution i.e. at some point the quality holds. Thus, solving $\psi_{\lambda}(w)=\lvert u \rvert$ is equivalent to
	\begin{align}
	2w^{3} + (4-2\lvert u \rvert) w^{2}+(2-4\lvert u \rvert)w + \lambda - 2\lvert u \rvert = 0.
	\label{eq:proxrat1}
	\end{align}
	Equation \eqref{eq:proxrat1} is easily solved using traditional python packages\footnote{Example of Python function to solve Eq. \eqref{eq:proxrat1} at \url{https://numpy.org/doc/stable/reference/generated/numpy.roots.html}. }.
	
	\begin{table}[t]
		\centering
		\caption{Invex regularization functions from Table \ref{tab:list} and their corresponding proximity operator ($\lambda \in (0,1]$ is a thresholding parameter).}
		\begin{tabular}{|p{0.5cm}| p{5.2cm}| p{6.7cm}|}
			\hline
			Ref & Invex function& Proximal operator \\
			\hline
			\cite{marjanovic2012l_q} & $g_{\lambda}(x)=\lambda \lvert x \rvert^{p}$, $p\in (0,1)$, $x\not =0$. & $\text{Prox}_{g_{\lambda}}(t)=\left\lbrace \begin{array}{ll}
			0 & \lvert t \rvert < \tau \\
			\{0,\text{sign}(t)\beta\} & \lvert t \rvert=\tau \\
			\text{sign}(t)y & \lvert t \rvert>\tau
			\end{array} \right.$ \vspace{0.5em}
			
			where $\beta = [2\lambda (1-p)]^{1/(2-p)}$, $\tau = \beta + \lambda p \beta^{p-1}$, $h(y)= \lambda p y^{p-1} + y - \lvert t \rvert = 0$, $y\in [\beta,\lvert t \rvert]$\\
			\hline
			
			- & $g_{\lambda}(x)= \lambda\log(1 + \lvert x \rvert)$ & $\text{Prox}_{g_{\lambda}}(t)=\left\lbrace \begin{array}{ll}
			0 & (\lvert t \rvert+1)^2< 4\lambda \\
			\text{sign}(t)\beta & \beta\geq 0 \\
			0 & \text{ otherwise }
			
		\end{array} \right.$ \vspace{0.5em}
		
		where $\beta = \frac{\lvert t \rvert -1 + \sqrt{(\lvert t \rvert + 1)^2 - 4\lambda}}{2}$.\\
		\hline
		
		- & $g_{\lambda}(x)=\lambda \frac{\lvert x \rvert}{2 + 2\lvert x \rvert}$ & $\text{Prox}_{g_{\lambda}}(t)=\left\lbrace \begin{array}{ll}
		0 & \lvert t \rvert = 0 \\
		\text{sign}(t)\beta & \text{ otherwise }
		\end{array} \right.$ \vspace{0.5em}
		
		where $2\beta^{3} + (4-2\lvert t \rvert) \beta^{2}+(2-4\lvert t \rvert)\beta + \lambda - 2\lvert t \rvert = 0$, $\beta>0$, and closest to $\lvert t \rvert$.\\
		\hline
		
		- & $g_{\lambda}(x)=\lambda \frac{x^{2}}{1 + x^{2}}$ & $\text{Prox}_{g_{\lambda}}(t)=\left\lbrace \begin{array}{ll}
		0 & \lvert t \rvert = 0 \\
		\text{sign}(t)\beta & \text{ otherwise }
		\end{array} \right.$ \vspace{0.5em}
		
		where $\beta^{5} - \lvert t \rvert \beta^{4} + 2\beta^{3} - 2\lvert t \rvert \beta^{2}+(1+2\lambda)\beta - \lvert t \rvert = 0$, $\beta>0$, and closest to $\lvert t \rvert$\\
		\hline
		
		- & $g_{\lambda}(x)=\lambda \left(\log(1 + \lvert x \rvert) - \frac{\lvert x \rvert}{2 + 2\lvert x \rvert}\right)$ & $\text{Prox}_{g_{\lambda}}(t)=\left\lbrace \begin{array}{ll}
		0 & \lvert t \rvert = 0 \\
		\text{sign}(t)\beta & \text{ otherwise } 
		\end{array} \right.$ \vspace{0.5em}
		
		where $2\beta^{3} + (4-2\lvert t \rvert)\beta^{2} + (2\lambda + 2 -4\lvert t \rvert)\beta + \lambda - 2\lvert t \rvert=0$, $\beta >0$, and closest to $\lvert t \rvert$.\\
		\hline
	\end{tabular}
	\label{tab:proximals}
\end{table}

\paragraph{Proximal of Eq. \eqref{fun4}}
Consider $h(w)=\lambda \frac{w^{2}}{1 + w^{2}} + \frac{1}{2}(w-u)^{2}$ for $\lambda\in (0,1]$, and fixed $u\in \mathbb{R}$. We note first that we only consider $w's$ for which $\text{sign}(w)=\text{sign}(u)$, otherwise $h(w)=\lambda \frac{w^{2}}{1 + w^{2}} + \frac{1}{2}w^{2}+\lvert u \rvert\lvert w \rvert + \frac{1}{2}u^{2}$ which is clearly minimized at $w=0$. Then, since with $\text{sign}(w)=\text{sign}(u)$ we have $(w-u)^{2} = (\lvert w \rvert-\lvert u \rvert)^{2}$, we replace $u$ with $\lvert u \rvert$ and take $w\geq 0$. As $h(w)$ is differentiable for $w>0$, re-arranging $h^{\prime}(w)=0$ gives
\begin{align}
\psi_{\lambda}(w)\delequal\frac{2\lambda w}{(1+w^{2})^{2}} + w = \lvert u \rvert.
\label{eq:proxrat2}
\end{align}
Observe that $\psi^{\prime}_{\lambda}(w)$ is always positive then it means that $\psi_{\lambda}(w)$ is monotonically increasing. Thus, the equation $\psi_{\lambda}(w)=\lvert u \rvert $ has unique solution i.e. at some point the quality holds. Thus, solving $\psi_{\lambda}(w)=\lvert u \rvert$ is equivalent to
\begin{align}
w^{5} - \lvert u \rvert w^{4} + 2w^{3} - 2\lvert u \rvert w^{2}+(1+2\lambda)w - \lvert u \rvert = 0.
\label{eq:proxrat3}
\end{align}
Equation \eqref{eq:proxrat3} is easily solved using traditional python packages.

\paragraph{Proximal of Eq. \eqref{fun5}}
Consider $h(w)=\lambda \left(\log(1 + \lvert w \rvert) - \frac{\lvert w \rvert}{2 + 2\lvert w \rvert}\right) + \frac{1}{2}(w-u)^{2}$ for $\lambda\in (0,1]$, and fixed $u\in \mathbb{R}$. We note first that we only consider $w's$ for which $\text{sign}(w)=\text{sign}(u)$, otherwise $h(w)=\lambda \left(\log(1 + \lvert w \rvert) - \frac{\lvert w \rvert}{2 + 2\lvert w \rvert}\right) + \frac{1}{2}w^{2}+\lvert u \rvert\lvert w \rvert + \frac{1}{2}u^{2}$ which is clearly minimized at $w=0$. Then, since with $\text{sign}(w)=\text{sign}(u)$ we have $(w-u)^{2} = (\lvert w \rvert-\lvert u \rvert)^{2}$, we replace $u$ with $\lvert u \rvert$ and take $w\geq 0$. As $h(w)$ is differentiable for $w>0$, re-arranging $h^{\prime}(w)=0$ gives
\begin{align}
\psi_{\lambda}(w)\delequal\lambda\frac{2w+1}{2(1+w)^{2}} + w = \lvert u \rvert.
\label{eq:proxlog2}
\end{align}
Observe that $\psi^{\prime}_{\lambda}(w)$ is always positive then it means that $\psi_{\lambda}(w)$ is monotonically increasing. Thus, the equation $\psi_{\lambda}(w)=\lvert u \rvert $ has unique solution i.e. at some point the quality holds. Thus, solving $\psi_{\lambda}(w)=\lvert u \rvert$ is equivalent to
\begin{align}
2w^{3} + (4-2\lvert u \rvert)w^{2} + (2\lambda + 2 -4\lvert u \rvert)w + \lambda - 2\lvert u \rvert=0.
\label{eq:proxlog3}
\end{align}
Equation \eqref{eq:proxlog3} is easily solved using traditional python packages.

\section{Proof of Theorem \ref{theo:ourCS}}
\label{app:ourCS}
We split the proof of Theorem \ref{theo:ourCS} into two parts. First part we focus our analysis on functions in Eqs. \eqref{fun2},\eqref{fun3},\eqref{fun5} and second part Eq. \eqref{fun4} (extra mild conditions are needed). Recall we skipped Eq. \eqref{fun1}.

\subsection{Part one}
We particularized \cite[Theorem 1]{gribonval2007highly} in order to prove that whenever the $\ell_{1}$-norm solution of optimization problem in Eq. \eqref{eq:problem4} is unique, then Eq. \eqref{eq:problem4} when $g(\boldsymbol{x})$ satisfies the following definition has the same global optima.

\begin{definition}{(\textit{Sparseness measure} \cite{gribonval2007highly})}
Let $g:\mathbb{R}^{n}\rightarrow \mathbb{R}$ such that $g(\boldsymbol{w})=\sum_{i=1}^{n}r(\boldsymbol{w}[i])$, where $r:[0,\infty)\rightarrow [0,\infty)$ and increasing. If $r$, not identically zero, with $r(0)=0$ such that $r(t)/t$ is non-increasing on $(0,\infty)$, then $g(\boldsymbol{x})$ is said to be a \textit{sparseness measure}.
\label{def:measure}
\end{definition}

Now we present the particular version in \cite[Theorem 1]{gribonval2007highly} as follows.

\begin{lemma}
Assume $\boldsymbol{H}\boldsymbol{x}=\boldsymbol{b}$, where $\boldsymbol{x}\in \mathbb{R}^{n}$ is $k$-sparse, the matrix $\boldsymbol{H}\in \mathbb{R}^{m\times n}$ ($m<n$) with $\ell_{2}$-normalized columns that satisfies RIP for any $2k$-sparse vector, with $\delta_{2k}<\frac{1}{3}$, and $\boldsymbol{b}\in \mathbb{R}^{m}$ is a noiseless measurements data vector. If $g(\boldsymbol{x})$ in Eqs. \eqref{fun2},\eqref{fun3}, \eqref{fun5} satisfies Definition \ref{def:measure}, then $\boldsymbol{x}$ is exactly recovered by solving Eq. \eqref{eq:problem4} i.e. only global optimizers exists.
\end{lemma}

In the following we prove functions in Eqs. \eqref{fun2},\eqref{fun3}, and \eqref{fun5} satisfy Definition \ref{def:measure}, and we proceed by cases.

\begin{wrapfigure}[14]{r}{0.4\textwidth}
	\vspace{-2em}
	\begin{center}
		\includegraphics[width=0.38\textwidth]{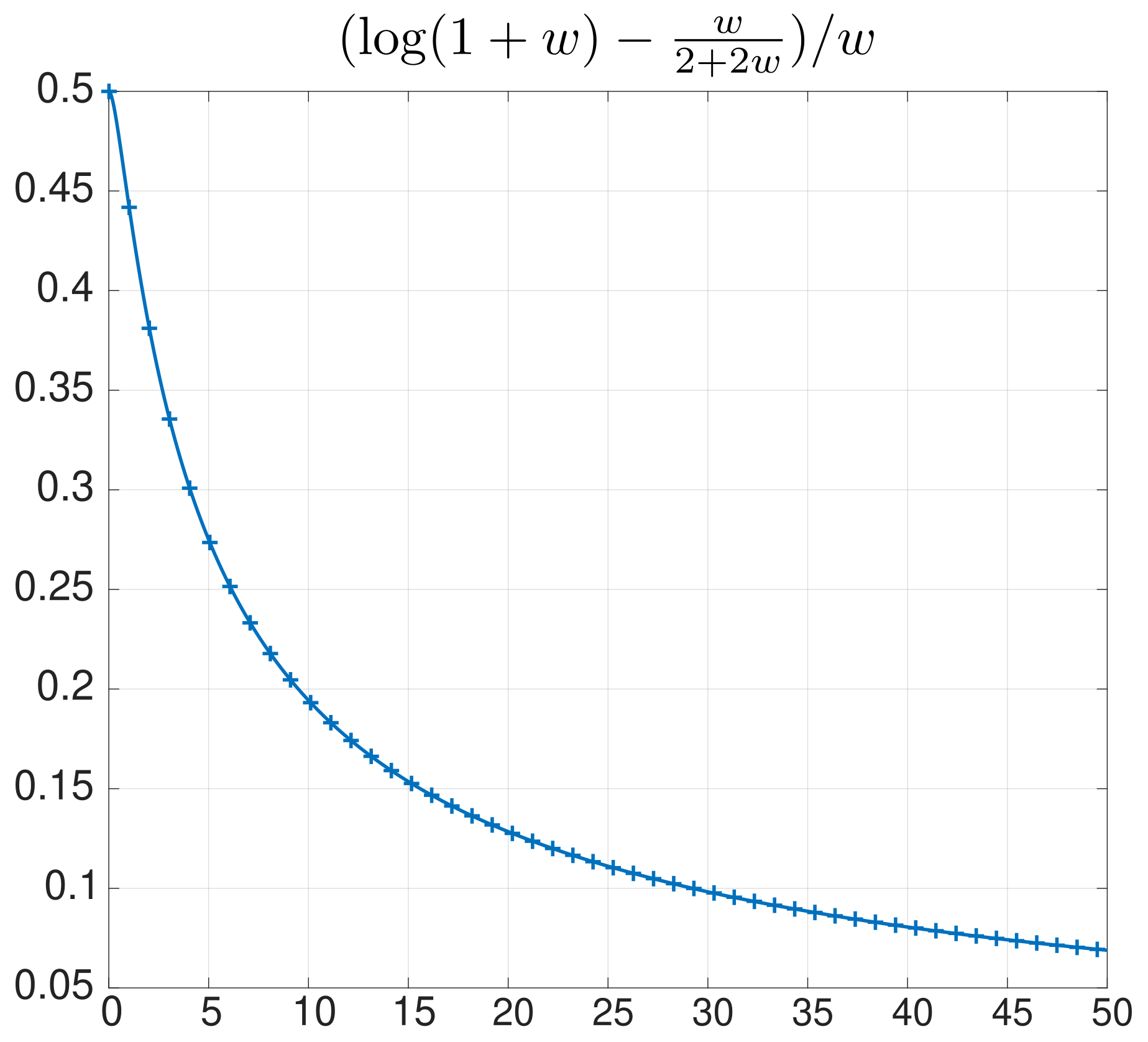}
	\end{center}
	\caption{\small Plot of $g(w)/w$ for $g(w)$ being Eq. \eqref{fun5} and $w>0$ to check that $g(w)/w$ is non-increasing on $(0,\infty)$.}
	\label{fig:auxLogProof}
\end{wrapfigure}

\begin{proof}
	
\textbf{Eq. \eqref{fun2}: }Take $g(w)= \log(1 + \lvert w \rvert)$ for any $w\in \mathbb{R}$. It is trivial to see that $g(0)=0$, and that $g(w)$ it is not identically zero. Then, we just need to show that $g(w)/w$ is non-increasing on $(0,\infty)$. Define $r(w)=\frac{\log(1+w)}{w}$. Observe that the derivative of $r(w)$ is given by $r^{\prime}(w) = \frac{\frac{w}{1+w}-\log(1+w)}{w^{2}}$, for $w\in (0,\infty)$. Since $\frac{w}{1+w}-\log(1+w)< 0$, then $r^{\prime}(w)< 0$ leads to conclude that $g(w)/w$ is non-increasing on $(0,\infty)$.

\textbf{Eq. \eqref{fun3}: }Take $g(w)= \frac{\lvert w\rvert}{2 + 2\lvert w\rvert}$ for any $w\in \mathbb{R}$. It is trivial to see that $g(0)=0$, and that $g(w)$ it is not identically zero. Then, we just need to show that $g(w)/w$ is non-increasing on $(0,\infty)$. Define $r(w)=\frac{w}{2w+2w^{2}} = \frac{1}{2+2w}$. Then, it is clear to conclude that $g(w)/w$ is non-increasing on $(0,\infty)$.

\textbf{Eq. \eqref{fun5}: }Take $g(w)= \log(1+\lvert w \rvert) - \frac{\lvert w \rvert}{2 + 2\lvert w \rvert}$ for any $w\in \mathbb{R}$. It is trivial to see that $g(0)=0$, and that $g(w)$ it is not identically zero. Then, we just need to show that $g(w)/w$ is non-increasing on $(0,\infty)$. For easy of exposition we present in Figure \ref{fig:auxLogProof} the plot of $g(w)/w$. Then it is clear that $g(w)/w$ is non-increasing on $(0,\infty)$.
\end{proof}

\subsection{Part two}
For this second part we appeal to a generalized result of \cite[Theorem 1]{gribonval2007highly} presented in \cite[Theorem 3.10]{woodworth2016compressed}. To exploit this generalized theorem we introduce the following definition.

\begin{definition}{(\textit{Admissible sparseness measure} \cite{woodworth2016compressed})}
A function $g:\mathbb{R}^{n}\rightarrow \mathbb{R}$ such that $g(\boldsymbol{w})=\sum_{i=1}^{n}r(\boldsymbol{w}[i])$ is said to be an admissible sparseness measure if 
\begin{itemize}
	\item $r(0)=0$, and g even on $\mathbb{R}$,
	\item $r$ is continuous on $\mathbb{R}$, and strictly increasing and strictly concave on $\mathbb{R}$.
\end{itemize}
\label{def:admissible}
\end{definition}

Based on the above definition we particularized \cite[Theorem 3.10]{woodworth2016compressed} in the lemma below in order to prove the solution of optimization problem in Eq. \eqref{eq:problem4} is unique, when functions in Eq. \eqref{fun4} are used under some mild conditions.

\begin{lemma}{(\cite[Theorem 3.10]{woodworth2016compressed})}
Assume $\boldsymbol{H}\boldsymbol{x}=\boldsymbol{b}$, where $\boldsymbol{x}\in \mathbb{R}^{n}$ is $k$-sparse, the matrix $\boldsymbol{H}\in \mathbb{R}^{m\times n}$ ($m<n$) with $\ell_{2}$-normalized columns that satisfies RIP for $\delta_{s}\in (0,1)$ with $s\geq 2k$, and $\boldsymbol{b}\in \mathbb{R}^{m}$ is a noiseless measurements data vector. Define $\beta_{1},\beta_{2} >0$ to be the lower and upper bound of magnitudes of non-zero entries of feasible vectors of Eq. \eqref{eq:problem4} (their existence if guaranteed \cite{woodworth2016compressed}). If $kr(2\beta_{2}) < (s+k-1)r(\beta_{1})$, then $\boldsymbol{x}$ is exactly recovered by solving Eq. \eqref{eq:problem4} i.e. only global optimizers exists.
\end{lemma}	

In the following we prove functions from Eq. \eqref{fun4} in Table \ref{tab:list} are able to exactly recover the signal $\boldsymbol{x}$ under some mild conditions.
\begin{proof}

\textbf{Eq. \eqref{fun4}: }Take $r(w)= \frac{w^{2}}{1 + w^{2}}$ for any $w\in \mathbb{R}$. It is trivial to see that $r(0)=0$, to check that it is even, continuous, and strictly increasing. Observe that the second derivative of $r(w)$ is given by $r^{\prime \prime}(w)=\frac{2-6w^{2}}{(1+w^2)^2}$. Then it is clear to conclude that $r(w)$ is strictly concave when $w>\frac{1}{3}$. Then, in order to have the chance to exactly recover the signal $\boldsymbol{x}$ we need to assume that the lower bound of magnitudes of non-zero entries of feasible vectors is $\beta_{1}>\frac{1}{3}$. Without loss of generality we assume $\boldsymbol{x}$ is a normalized signal (in practical imaging applications $\boldsymbol{x}$ is always normalized). Then, we take $\beta_{1} = 0.5$, and $\beta_{2}=1.0$. In addition, assuming $\boldsymbol{H}$ satisfies RIP when $s\geq 4k+2$, with $\delta_{s}\in (0,1)$, it is numerically easy to verified that $kr(2\beta_{2}) < (s+k-1)r(\beta_{1})$.
\end{proof}

\section{Proof of Lemma \ref{lem:convergeAPG}}
\label{app:lemAPG}
Before proving Lemma \ref{lem:convergeAPG} we consider two definitions in the following which the loss function $F(\boldsymbol{x})=f(\boldsymbol{x}) + \lambda g(\boldsymbol{x})$ in Eq. \eqref{eq:problem4} satisfies. Recall that $\lambda\in (0,1]$.

\begin{definition}
A function $h:\mathbb{R}^{n}\rightarrow (-\infty,\infty]$ is said to be proper if $\text{dom }h \not = \emptyset$, where $\text{dom }=\{\boldsymbol{x}\in \mathbb{R}^{n}: h(\boldsymbol{x})< \infty\}$.
\end{definition}
Since we are assuming the sensing matrix $\boldsymbol{H}$ satisfies RIP it guarantees the existence of a solution to Eq. \eqref{eq:problem4} implying that $\text{dom }F \not = \emptyset$. Thus, $F(\boldsymbol{x})$ in Eq. \eqref{eq:problem4} satisfies the above definition because.

\begin{definition}
A function $h:\mathbb{R}^{n}\rightarrow \mathbb{R}$ is coercive, if $h$ is bounded from below and $h(\boldsymbol{x})\rightarrow \infty$ when $\lVert \boldsymbol{x} \rVert_{2}\rightarrow \infty$.
\end{definition}
Considering that the list of invex functions in Table \ref{tab:list}, and $f(\boldsymbol{x})=\lVert \boldsymbol{H}\boldsymbol{x} -\boldsymbol{v} \rVert_{2}^{2}$ (for fix $\boldsymbol{v}$ and $\boldsymbol{H}$ satisfying RIP) are positive, then the loss function $F(\boldsymbol{x})$ satisfies $F(\boldsymbol{x})\geq 0$. The second part of the coercive definition is trivially guaranteed since $\boldsymbol{H}$ satisfies RIP, otherwise we will be denying the existence of a global solution to Eq. \eqref{eq:problem4} which is a contradiction.

Now we proceed to prove Lemma \ref{lem:convergeAPG}.

\begin{proof}
Line 6 in Algorithm \ref{alg:invexProximal} is given by
\begin{align}
\boldsymbol{v}^{(t+1)} = \argmin_{\boldsymbol{x}\in \mathbb{R}^{n}} \hspace{0.5em} \left\langle \nabla f(\boldsymbol{x}^{(t)}),\boldsymbol{x}-\boldsymbol{x}^{(t)} \right\rangle + \frac{1}{2\lambda \alpha_{1}} \lVert \boldsymbol{x}-\boldsymbol{x}^{(t)} \rVert_{2}^{2} + g(\boldsymbol{x}).
\label{eq:noninvex1}
\end{align}
We write equal in the above equation because the proximal in Line 6 is invex therefore it always map to a global optimizer. So from Eq. \eqref{eq:noninvex1} we have 
\begin{align}
\left\langle \nabla f(\boldsymbol{x}^{(t)}),\boldsymbol{v}^{(t+1)}-\boldsymbol{x}^{(t)} \right\rangle + \frac{1}{2\lambda\alpha_{1}} \lVert \boldsymbol{v}^{(t+1)}-\boldsymbol{x}^{(t)} \rVert_{2}^{2} + g(\boldsymbol{v}^{(t+1)}) \leq g(\boldsymbol{x}^{(t)}).
\label{eq:noninvex2}
\end{align}
From the Lipschitz continuous of $\nabla f$ and Eq. \eqref{eq:noninvex2} we have
\begin{align}
F(\boldsymbol{v}^{(t+1)}) &\leq g(\boldsymbol{v}^{(t+1)}) + f(\boldsymbol{x}^{(t)}) + \left\langle \nabla f(\boldsymbol{x}^{(t)}),\boldsymbol{v}^{(t+1)}-\boldsymbol{x}^{(t)} \right\rangle + \frac{L}{2} \lVert \boldsymbol{v}^{(t+1)}-\boldsymbol{x}^{(t)} \rVert_{2}^{2} \nonumber\\
&\leq g(\boldsymbol{x}^{(t)}) - \left\langle \nabla f(\boldsymbol{x}^{(t)}),\boldsymbol{v}^{(t+1)}-\boldsymbol{x}^{(t)} \right\rangle - \frac{1}{2\lambda\alpha_{1}} \lVert \boldsymbol{v}^{(t+1)}-\boldsymbol{x}^{(t)} \rVert_{2}^{2} \nonumber\\
&+ f(\boldsymbol{x}^{(t)}) + \left\langle \nabla f(\boldsymbol{x}^{(t)}),\boldsymbol{v}^{(t+1)}-\boldsymbol{x}^{(t)} \right\rangle + \frac{L}{2} \lVert \boldsymbol{v}^{(t+1)}-\boldsymbol{x}^{(t)} \rVert_{2}^{2} \nonumber\\
&=F(\boldsymbol{x}^{(t)}) - \left(\frac{1}{2\lambda\alpha_{1}}-\frac{L}{2}\right)\lVert \boldsymbol{v}^{(t+1)}-\boldsymbol{x}^{(t)} \rVert_{2}^{2}.
\label{eq:noninvex3}
\end{align}
If $F(\boldsymbol{z}^{(t+1)})\leq F(\boldsymbol{v}^{(t+1)})$, then
\begin{align}
\boldsymbol{x}^{(t+1)} = \boldsymbol{z}^{(t+1)}, F(\boldsymbol{x}^{(t+1)})=F(\boldsymbol{z}^{(t+1)})\leq F(\boldsymbol{v}^{(t+1)}).
\label{eq:noninvex4}
\end{align}
If $F(\boldsymbol{z}^{(t+1)})> F(\boldsymbol{v}^{(t+1)})$, then
\begin{align}
\boldsymbol{x}^{(t+1)} = \boldsymbol{v}^{(t+1)}, F(\boldsymbol{x}^{(t+1)})=F(\boldsymbol{v}^{(t+1)}).
\label{eq:noninvex5}
\end{align}
From Eqs. \eqref{eq:noninvex3}, \eqref{eq:noninvex4} and \eqref{eq:noninvex5} we have
\begin{align}
F(\boldsymbol{x}^{(t+1)} )\leq F(\boldsymbol{v}^{(t+1)}) \leq F(\boldsymbol{x}^{(t)}).
\label{eq:nonIncreasing}
\end{align}
So
\begin{align}
F(\boldsymbol{x}^{(t+1)})\leq F(\boldsymbol{x}^{(1)}), F(\boldsymbol{v}^{(t+1)})\leq F(\boldsymbol{x}^{(1)}),
\label{eq:noninvex6}
\end{align}
for all $t$. Recall that we consider the estimation of $\boldsymbol{z}^{(t+1)}$ unique because it is performed through the proximal of $g(\boldsymbol{x})$ which always map to a global optimizer.

Observe that from Eq. \eqref{eq:nonIncreasing} was concluded that $F(\boldsymbol{x}^{(t)})$ is nonincreasing then for all $t>1$ we have $F(\boldsymbol{x}^{(t)})\leq F(\boldsymbol{x}^{(1)})$ and therefore $\boldsymbol{x}^{(t)}\in \{\boldsymbol{w}: F(\boldsymbol{w})\leq F(\boldsymbol{x}^{(1)})\}$ (known as level sets). Since $F(\boldsymbol{x})$ is coercive then all its level sets are bounded. Then we know that $\{\boldsymbol{x}^{(t)} \}$, and $\{\boldsymbol{v}^{(t)}\}$ are also bounded. Thus $\{\boldsymbol{x}^{(t)} \}$ has accumulation points. Let $\boldsymbol{x}^{*}$ be any accumulation point of $\{\boldsymbol{x}^{(t)}\}$, say a subsequence satisfying $\{\boldsymbol{x}^{(t_{j}+1)}\}\rightarrow \boldsymbol{x}^{*}$ as $j\rightarrow \infty$. Let $F^{*}$ be $\displaystyle \lim_{j\rightarrow \infty} \hspace{0.5em}F(\boldsymbol{x}^{(t_{j}+1)}) = F(\boldsymbol{x}^{*})=F^{*}$. The existence of this limit is guaranteed since $f$ is continuously differentiable. Then, from Eq. \eqref{eq:noninvex3} we have
\begin{align}
\left(\frac{1}{2\lambda\alpha_{1}}-\frac{L}{2}\right)\lVert \boldsymbol{v}^{(t+1)}-\boldsymbol{x}^{(t)} \rVert_{2}^{2} \leq F(\boldsymbol{x}^{(t)}) - F(\boldsymbol{v}^{(t+1)}) \leq F(\boldsymbol{x}^{(t)}) - F(\boldsymbol{x}^{(t+1)}).
\label{eq:noninvex7}
\end{align}
Summing over $t=1,2,\dots,\infty$, we have
\begin{align}
\left(\frac{1}{2\lambda\alpha_{1}}-\frac{L}{2}\right)\sum_{t=1}^{\infty}\lVert \boldsymbol{v}^{(t+1)}-\boldsymbol{x}^{(t)} \rVert_{2}^{2} \leq F(\boldsymbol{x}^{(1)}) - F^{*} < \infty.
\label{eq:noninvex8}
\end{align}
From $\alpha_{1}< \frac{1}{L}$ we have
\begin{align}
\lVert \boldsymbol{v}^{(t+1)}-\boldsymbol{x}^{(t)} \rVert_{2}^{2}\rightarrow 0, \text{ as } t\rightarrow \infty.
\label{eq:noninvex9}
\end{align}
From the optimality condition of Eq. \eqref{eq:noninvex1} we have
\begin{align}
\boldsymbol{0} &\in \nabla f(\boldsymbol{x}^{(t)}) + \frac{1}{\lambda\alpha_{1}}(\boldsymbol{v}^{(t+1)} - \boldsymbol{x}^{(t)}) + \partial g(\boldsymbol{v}^{(t+1)}) \nonumber\\
&=\nabla f(\boldsymbol{x}^{(t)}) + \nabla f(\boldsymbol{v}^{(t+1)}) - \nabla f(\boldsymbol{v}^{(t+1)}) + \frac{1}{\lambda\alpha_{1}}(\boldsymbol{v}^{(t+1)} - \boldsymbol{x}^{(t)}) + \partial g(\boldsymbol{v}^{(t+1)}).
\label{eq:noninvex10}
\end{align}
So we have
\begin{align}
-\nabla f(\boldsymbol{x}^{(t)}) + \nabla f(\boldsymbol{v}^{(t+1)}) - \frac{1}{\lambda\alpha_{1}}(\boldsymbol{v}^{(t+1)} - \boldsymbol{x}^{(t)}) \in \partial F(\boldsymbol{v}^{(t+1)}),
\label{eq:noninvex11}
\end{align}
and
\begin{align}
\left \lVert \nabla f(\boldsymbol{x}^{(t)})- \nabla f(\boldsymbol{v}^{(t+1)}) + \frac{1}{\lambda\alpha_{1}}(\boldsymbol{v}^{(t+1)} - \boldsymbol{x}^{(t)}) \right \rVert_{2} \leq \left(\frac{1}{\lambda\alpha_{1}} + L\right) \lVert \boldsymbol{v}^{(t+1)} - \boldsymbol{x}^{(t)} \rVert_{2} \rightarrow 0,
\label{eq:noninvex12}
\end{align}
as $t\rightarrow \infty$.

From Eq. \eqref{eq:noninvex9} we have $\boldsymbol{v}^{(t_{j}+1)}\rightarrow \boldsymbol{x}^{*}$ as $j\rightarrow \infty$. From Eq. \eqref{eq:noninvex1} we have
\begin{align}
&\left\langle \nabla f(\boldsymbol{x}^{(t_{j})}),\boldsymbol{v}^{(t_{j}+1)}-\boldsymbol{x}^{(t_{j}+1)} \right\rangle + \frac{1}{2\lambda\alpha_{1}} \lVert \boldsymbol{v}^{(t_{j}+1)}-\boldsymbol{x}^{(t_{j})} \rVert_{2}^{2} + g(\boldsymbol{v}^{(t_{j}+1)}) \nonumber\\
&\leq \left\langle \nabla f(\boldsymbol{x}^{(t_{j})}),\boldsymbol{x}^{*}-\boldsymbol{x}^{(t_{j})} \right\rangle + \frac{1}{2\lambda\alpha_{1}} \lVert \boldsymbol{x}^{*}-\boldsymbol{x}^{(t_{j})} \rVert_{2}^{2} + g(\boldsymbol{x}^{*})
\label{eq:noninvex13}
\end{align}		
So
\begin{align}
\limsup_{j\rightarrow \infty} \hspace{0.5em}g(\boldsymbol{v}^{(t_{j}+1)}) \leq g(\boldsymbol{x}^{*}).
\label{eq:noninvex14}
\end{align}
From the continuity assumption on $g$ we have $\displaystyle \liminf_{j\rightarrow \infty} \hspace{0.5em}g(\boldsymbol{v}^{(t_{j}+1)}) \geq g(\boldsymbol{x}^{*})$, then we conclude
\begin{align}
\lim_{j\rightarrow \infty} \hspace{0.5em}g(\boldsymbol{v}^{(t_{j}+1)}) = g(\boldsymbol{x}^{*}).
\label{eq:noninvex15}
\end{align}
Because $f$ is continuously differentiable, we have $\displaystyle \lim_{j\rightarrow \infty} \hspace{0.5em}F(\boldsymbol{v}^{(t_{j}+1)}) = F(\boldsymbol{x}^{*})$. From $\{\boldsymbol{v}^{(t_{j}+1)}\}\rightarrow \boldsymbol{x}^{*}$, and Eq. \eqref{eq:noninvex11} we have $\boldsymbol{0} \in \partial F(\boldsymbol{x}^{*})$. Therefore, since $F(\boldsymbol{x})$ is invex according to Theorem \ref{theo:ourCS} we have that the sequence $\{\boldsymbol{x}^{(t)}\}$ converges to a global minimizer of $F(\boldsymbol{x})$.
\end{proof}

\subsection{Numerical Validation of Lemma \ref{lem:convergeAPG}}
\begin{figure}[h!]
	\centering
	\includegraphics[width=1\linewidth]{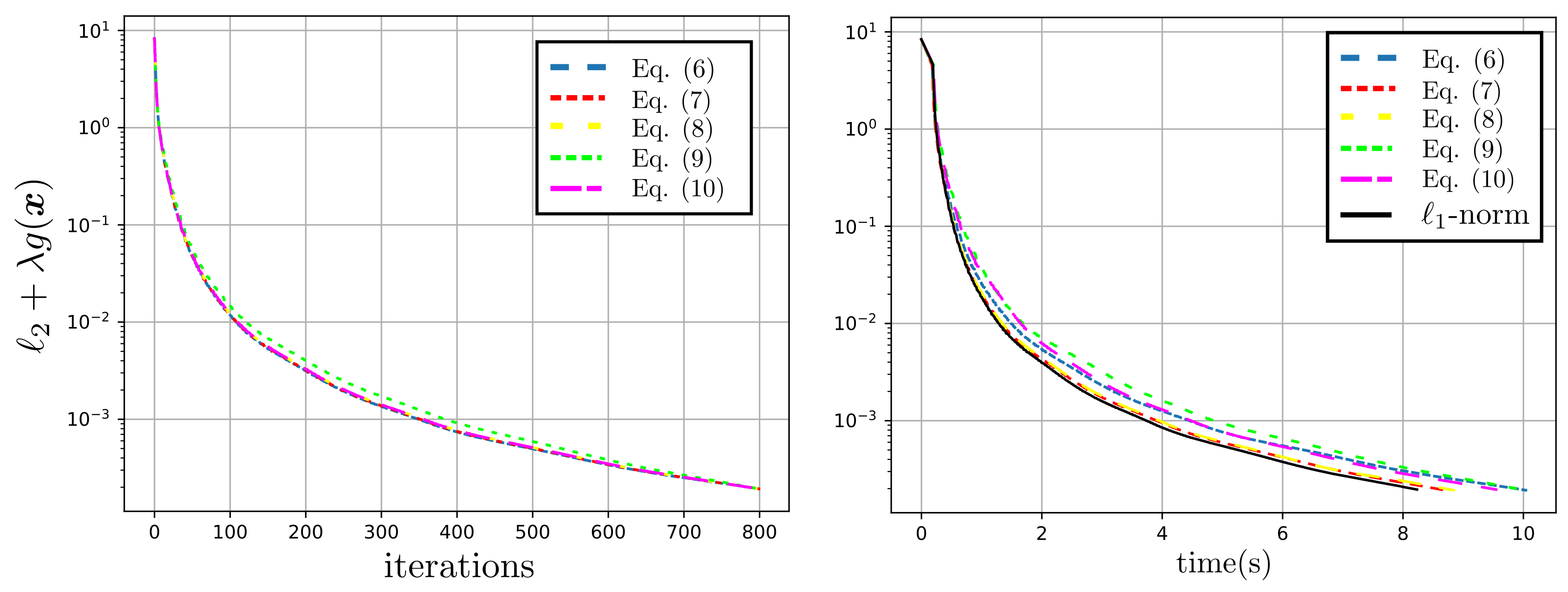}
	\vspace{-0.3em}
	\caption{Numerical convergence of Algorithm \ref{alg:invexProximal} for all invex regularizers to recover an image of size $256\times 256$ from blurred data, as explained in Experiment 1. (left) Minimization process of $F(\boldsymbol{x})=\ell_{2}+ \lambda g(\boldsymbol{x})$ along $T=800$ iterations. (right) running time of Algorithm \ref{alg:invexProximal} to perform $T=800$ iterations for all invex regularizers and the $\ell_{1}$-norm.}
	\label{fig:msevstimeiter}
\end{figure}

To numerically validate the proof of Lemma \ref{lem:convergeAPG} provided in the above section, we present Fig. \ref{fig:msevstimeiter}. In this figure we are reporting the numerical convergence of Algorithm \ref{alg:invexProximal} for all invex regularizers to recover an image of size $256\times 256$ from blurred data, as explained in Experiment 1 for noiseless case. Specifically, Fig. \ref{fig:msevstimeiter}(left) reports how the loss function $F(\boldsymbol{x})=\ell_{2}+ \lambda g(\boldsymbol{x})$, analyzed in the above proof, is minimized along $T=800$ iterations. This plot numerically validates the proof of Lemma \ref{lem:convergeAPG}. As a complement to this plot, Fig. \ref{fig:msevstimeiter}(right) presents the running time of Algorithm \ref{alg:invexProximal} to perform $T=800$ iterations for all invex regularizers and the $\ell_{1}$-norm. This second plot suggests that Algorithm \ref{alg:invexProximal} using the $\ell_{1}$-norm as regularizer requires 1.8 seconds less than its invex competitors to perform $T=800$ iterations. We remark that this negligible difference is expected, since in Table \ref{tab:timesProx} was concluded that the running time to compute the proximal operator for all invex differs in the order of milliseconds with the computation of the proximal of $\ell_{1}$-norm.

\newpage
\section{Proof of Lemma \ref{theo:PnP}}
\label{app:PnP}
We proceed to prove this lemma by extending the mathematical analysis in \cite{kamilov2017plug} to invex functions. To that end, we recall some definitions and a classical result from monotone operator theory needed for the proof of this lemma as follows.

\begin{definition}{(\textit{Nonexpansiveness})}
An operator $F:\mathbb{R}^{n}\rightarrow \mathbb{R}^{n}$ is said to be nonexpansive if it is Lipschitz continuous as in Definition \ref{def:lipschitz} with $L=1$.
\label{def:nonexpansive}
\end{definition}	
Based on the nonexpansiveness concept we give the following definition.
\begin{definition}
For a constant $\beta \in (0,1)$ we say a function $G$ is $\beta$-average, if there exists a nonexpansive operator $F$ such that $G = (1-\beta)\boldsymbol{I} + \beta F$ 
\label{def:average}
\end{definition}
Now based on the concept of average operators we recall the following classical results.

\begin{lemma}{(\cite[Proposition 4.44]{bauschke2011convex})}
Let $G_{1}$ be $\beta_{1}$-averaged and $G_{2}$ be $\beta_{2}$-averaged. Then, the composite operator $G\delequal G_{2}\circ G_{1}$ is
\begin{align}
\beta \delequal \frac{\beta_{1} + \beta_{2} - 2\beta_{1}\beta_{2}}{1-\beta_{1}\beta_{2}},
\end{align}
averaged operator.
\label{lem:composedAverage}
\end{lemma}	

\begin{lemma}
Let $F$ be a $\beta$-average operator with $\beta \in (0,1)$. Then
\begin{align}
\lVert F(\boldsymbol{x})-F(\boldsymbol{y}) \rVert_{2}^{2} \leq \lVert \boldsymbol{x} -\boldsymbol{y} \rVert_{2}^{2} - \left(\frac{1-\beta}{\beta}\right) \lVert \boldsymbol{x}-F(\boldsymbol{x}) - \boldsymbol{y} + F(\boldsymbol{y})\rVert_{2}^{2}.
\end{align}
\label{lem:equiavalentAverage}
\end{lemma}

Now we proceed to prove Lemma \ref{theo:PnP}
\begin{proof}
Following the assumptions made in Lemma \ref{theo:PnP}, we start this proof by noticing that for a differentiable invex function $f$ a point $\boldsymbol{x}$ is a global minimizer of $f$ according to Theorem \ref{theo:optimal_v0} if 
\begin{align}
\boldsymbol{0} = \nabla f(\boldsymbol{x}) \leftrightarrow \boldsymbol{x} = (\boldsymbol{I}-\alpha \nabla f )(\boldsymbol{x}),
\end{align}
for non-zero $\alpha$. In other words, $\boldsymbol{x}$ is a minimizer of $f$ if and only if it is a fixed point of the mapping $\boldsymbol{I}-\alpha \nabla f$. This property of invex functions is what allows to extend the mathematical guarantees in \cite{kamilov2017plug} given only for convex functions. Now considering that $f$ is assumed to have Lipschitz continuous gradient with parameter $L$, then the operator $\boldsymbol{I}-\alpha \nabla f$ is Lipschitz with parameter $L_{G}=\max\{1,\lvert 1-\alpha L \rvert \}$ and therefore is nonexpansive for $\alpha \in (0,2/L]$. So it is averaged for $\alpha \in (0,2/L)$ since 
\begin{align}
\boldsymbol{I}-\alpha \nabla f = (1-\kappa) \boldsymbol{I} + \kappa \left(\boldsymbol{I}-2/L \nabla f\right),
\end{align}
where $\kappa = \alpha L/2 <1$.

Assume the denoiser $d$ is $\kappa$-averaged and the operator $G_{\alpha} = \boldsymbol{I}-\alpha \nabla f$. Observe that $G_{\alpha}$ is $(\gamma L/2)$-averaged for any $\alpha \in (0,2/L)$. From Lemma \ref{lem:composedAverage} their composition $P=d\circ G_{\alpha}$ is
\begin{align}
\beta \delequal \frac{\kappa + \gamma L/2 - 2\kappa\gamma L/2}{1-\kappa\gamma L/2},
\end{align}
averaged. Consider a single iteration $\boldsymbol{v}^{+}=P(\boldsymbol{x})$, then we have for any $\boldsymbol{x}^{*}=P(\boldsymbol{x}^{*})$ (fixed point) we have that
\begin{align}
\lVert \boldsymbol{v}^{+} - \boldsymbol{x}^{*} \rVert_{2}^{2} &¨= \lVert P(\boldsymbol{x}) - P(\boldsymbol{x}^{*}) \rVert_{2}^{2} \nonumber\\
&\leq \lVert \boldsymbol{x} - \boldsymbol{x}^{*} \rVert_{2}^{2} - \left(\frac{1-\beta}{\beta}\right) \lVert \boldsymbol{x}-P(\boldsymbol{x}) - \boldsymbol{x}^{*} + P(\boldsymbol{x}^{*}) \rVert_{2}^{2} \nonumber\\
&= \lVert \boldsymbol{x} - \boldsymbol{x}^{*} \rVert_{2}^{2} - \left(\frac{1-\beta}{\beta}\right) \lVert \boldsymbol{x}-P(\boldsymbol{x}) \rVert_{2}^{2},
\end{align}
where we used Lemma \ref{lem:equiavalentAverage}. From Line 6 in Algorithm \ref{alg:invexPnP} the iteration $t+1$ and rearranging the terms, we obtain
\begin{align}
\lVert \boldsymbol{x}^{(t)} - P(\boldsymbol{x}^{(t)})\rVert_{2}^{2} \leq \left(\frac{\beta}{1-\beta}\right) \left[\lVert \boldsymbol{x}^{(t)} - \boldsymbol{x}^{*} \rVert_{2}^{2} - \lVert \boldsymbol{v}^{(t+1)} - \boldsymbol{x}^{*} \rVert_{2}^{2} \right].
\end{align}
If $F(\boldsymbol{z}^{(t+1)})\leq F(\boldsymbol{v}^{(t+1)})$, then
\begin{align}
\boldsymbol{x}^{(t+1)} = \boldsymbol{z}^{(t+1)}, F(\boldsymbol{x}^{(t+1)})=F(\boldsymbol{z}^{(t+1)})\leq F(\boldsymbol{v}^{(t+1)}).
\label{eq:noninvex41}
\end{align}
If $F(\boldsymbol{z}^{(t+1)})> F(\boldsymbol{v}^{(t+1)})$, then
\begin{align}
\boldsymbol{x}^{(t+1)} = \boldsymbol{v}^{(t+1)}, F(\boldsymbol{x}^{(t+1)})=F(\boldsymbol{v}^{(t+1)}).
\label{eq:noninvex51}
\end{align}
From Eqs. \eqref{eq:noninvex41} and \eqref{eq:noninvex51} we have
\begin{align}
F(\boldsymbol{x}^{(t+1)} )\leq F(\boldsymbol{v}^{(t+1)}) \leq F(\boldsymbol{x}^{(t)}).
\label{eq:nonIncreasing1}
\end{align}

Observe that from Eq. \eqref{eq:nonIncreasing1} was concluded that $F(\boldsymbol{x}^{(t)})$ is nonincreasing then for all $t>1$ we have $F(\boldsymbol{x}^{(t)})\leq F(\boldsymbol{x}^{(1)})$ and therefore $\boldsymbol{x}^{(t)}\in \{\boldsymbol{w}: F(\boldsymbol{w})\leq F(\boldsymbol{x}^{(1)})\}$ (known as level sets). Since $F(\boldsymbol{x})$ is coercive then all its level sets are bounded (concluded from Appendix \ref{app:lemAPG}). Then we know that $\{\boldsymbol{x}^{(t)} \}$, and $\{\boldsymbol{v}^{(t)}\}$ are also bounded. Thus $\{\boldsymbol{x}^{(t)} \}$ has accumulation points which guarantees the existence of $\boldsymbol{x}^{*}$ implying that for a subsequence satisfying $\{\boldsymbol{x}^{(t_{j}+1)}\}\rightarrow \boldsymbol{x}^{*}$ as $j\rightarrow \infty$, we also have $\displaystyle \lim_{j\rightarrow \infty} \hspace{0.5em}F(\boldsymbol{x}^{(t_{j}+1)}) = F(\boldsymbol{x}^{*})=F^{*}$. Then, from Eqs. \eqref{eq:noninvex41}, \eqref{eq:noninvex51}, and the continuity of $F$, it is easy to see that $\lVert \boldsymbol{x}^{(t+1)} - \boldsymbol{x}^{*} \rVert_{2}^{2}\leq \lVert \boldsymbol{v}^{(t+1)} - \boldsymbol{x}^{*} \rVert_{2}^{2}$, which leads to 
\begin{align}
\lVert \boldsymbol{x}^{(t)} - P(\boldsymbol{x}^{(t)})\rVert_{2}^{2} \leq \left(\frac{\beta}{1-\beta}\right) \left[\lVert \boldsymbol{x}^{(t)} - \boldsymbol{x}^{*} \rVert_{2}^{2} - \lVert \boldsymbol{x}^{(t+1)} - \boldsymbol{x}^{*} \rVert_{2}^{2} \right].
\end{align}

By averaging this inequality over $T$ iterations and dropping the last term $\lVert \boldsymbol{x}^{(t+1)} - \boldsymbol{x}^{*} \rVert_{2}^{2}$, we obtain
\begin{align}
\frac{1}{T}\sum_{t=1}^{T} \lVert \boldsymbol{x}^{(t)} - P(\boldsymbol{x}^{(t)}) \rVert_{2}^{2} \leq \frac{2}{T}\left(\frac{1+\kappa}{1-\kappa}\right) \lVert \boldsymbol{x}^{(0)} - \boldsymbol{x}^{*} \rVert_{2}^{2}.
\label{eq:finalPnP}
\end{align}
To obtain the result that depends on $\kappa \in (0,1)$, we note that for any $\alpha \in (0,1/L]$, we write
\begin{align}
\frac{\beta}{1-\beta} = \frac{\kappa + \alpha L/2 - \kappa \alpha L}{(1-\kappa)(1-\alpha L/2)} \leq \frac{\kappa + \frac{1}{2}}{\frac{1-\kappa }{2}} \leq 2 \left(\frac{1+\kappa}{1-\kappa}\right).
\label{eq:finalPnP1}
\end{align}
Thus, from Eqs. \eqref{eq:finalPnP} and \eqref{eq:finalPnP1} the result holds.
\end{proof}

\subsection{Pseudo-code for plug-and-play invex imaging}
For the sake of completeness we present Algorithm \ref{alg:invexPnP} which is the pseudo-code of the plug-and-play version of APG for solving Eq. \eqref{eq:problem4}. The scaled-up convergence of APG are offered by two auxiliary variables, i.e., $\boldsymbol{y}^{(t+1)}$ and $\boldsymbol{z}^{(t+1)}$ in Lines 4 and 5. In Line 6 is presented the replacement of the proximal operator in APG pseudo-code with a neural network based denoiser Noise2Void \cite{krull2019noise2void}. And a monitor constrain computed in Line 8, to satisfy the sufficient descent property. 

\begin{algorithm}[H]
	\caption{Plug-and-play Proximal Gradient Algorithm}
	\label{alg:invexPnP}
	\begin{algorithmic}[1]
		\State{\textbf{input}: Tolerance constant $\epsilon\in (0,1)$, initial point $\boldsymbol{x}^{(0)}$, and number of iterations $T$}
		\State{\textbf{initialize}: $\boldsymbol{x}^{(1)}=\boldsymbol{x}^{(0)}=\boldsymbol{z}^{(0)}, r_{1}=1,r_{0}=0, \alpha_{1},\alpha_{2}< \frac{1}{L}$, and $\lambda \in (0,1]$}
		\For{$t=1$ to $T$}
		\State{$\boldsymbol{y}^{(t)}= \boldsymbol{x}^{(t)} + \frac{r_{t-1}}{r_{t}}(\boldsymbol{z}^{(t)}-\boldsymbol{x}^{(t)}) + \frac{r_{t-1}-1}{r_{t}}(\boldsymbol{x}^{(t)}- \boldsymbol{x}^{(t-1)})$}
		\State{$\boldsymbol{z}^{(t+1)}=\text{prox}_{\alpha_{2} \lambda g}(\boldsymbol{y}^{(t)} - \alpha_{2}\nabla f(\boldsymbol{y}^{(t)}))$}
		\State{$\boldsymbol{v}^{(t+1)} = \text{Noise2Void}(\boldsymbol{x}^{(t)} - \alpha_{1}\nabla f(\boldsymbol{x}^{(t)}))$ \Comment{This calls the trained Noise2Void model}}
		\State{$r_{t+1}=\frac{\sqrt{4(r_{t})^{2}+1}+1}{2}$}
		\State{$\boldsymbol{x}^{(t+1)}=\left \lbrace\begin{array}{ll}
				\boldsymbol{z}^{(t+1)}, & \text{ if }f(\boldsymbol{z}^{(t+1)} ) + \lambda g(\boldsymbol{z}^{(t+1)} )\leq f(\boldsymbol{v}^{(t+1)} ) + \lambda g(\boldsymbol{v}^{(t+1)} )\\
				\boldsymbol{v}^{(t+1)}, & \text{ otherwise }
			\end{array}\right.$}
		\EndFor
		\State{\textbf{return:} $\boldsymbol{x}^{(T)}$}
	\end{algorithmic}
\end{algorithm}

\section{Proof of Lemma \ref{lem:convergeUnrolling}}
\label{app:unrolling}
\begin{proof}
To prove this lemma we start exploiting the convexity of $\lambda g(\boldsymbol{x}) + \frac{1}{2}\lVert \boldsymbol{x} -\boldsymbol{u}\rVert^{2}_{2}$ for fixed $\boldsymbol{v}\in \mathbb{R}^{n}$ according to Theorem \ref{theo:proximalProof}, and $\lambda\in (0,1]$. Then, for all functions in Table \ref{tab:list}, we have
\begin{align}
\lambda g(\boldsymbol{x}) + \frac{1}{2}\lVert \boldsymbol{x} -\boldsymbol{u}\rVert^{2}_{2} - \lambda g(\boldsymbol{y}) - \frac{1}{2}\lVert \boldsymbol{y} -\boldsymbol{u}\rVert^{2}_{2} &\geq \left(\boldsymbol{\zeta}_{y} + \boldsymbol{y}-\boldsymbol{u}\right)^{T}(\boldsymbol{x}-\boldsymbol{y}) \nonumber\\
\lambda g(\boldsymbol{x}) - \lambda g(\boldsymbol{y}) &\geq \left(\boldsymbol{\zeta}_{y} + \boldsymbol{y}-\boldsymbol{u}\right)^{T}(\boldsymbol{x}-\boldsymbol{y}) + \frac{1}{2}\lVert \boldsymbol{y} -\boldsymbol{u}\rVert^{2}_{2} - \frac{1}{2}\lVert \boldsymbol{x} -\boldsymbol{u}\rVert^{2}_{2} \nonumber\\
\lambda g(\boldsymbol{x}) - \lambda g(\boldsymbol{y}) & \geq \left(\boldsymbol{\zeta}_{y} + \boldsymbol{y}-\boldsymbol{u}\right)^{T}(\boldsymbol{x}-\boldsymbol{y}) + (\boldsymbol{x}-\boldsymbol{u})^{T}(\boldsymbol{y}-\boldsymbol{x})
\label{eq:proofLem1}
\end{align}
for all $\boldsymbol{x},\boldsymbol{y}\in \mathbb{R}^{n}$, and $\boldsymbol{\zeta}_{y}\in \partial \lambda g(\boldsymbol{y})$, where the third inequality comes from the convexity of $f(\boldsymbol{x}) = \frac{1}{2}\lVert \boldsymbol{x} -\boldsymbol{u}\rVert^{2}_{2}$. Then, from Eq. \eqref{eq:proofLem1} we conclude
\begin{align}
\lambda g(\boldsymbol{x}) - \lambda g(\boldsymbol{y}) & \geq \boldsymbol{\zeta}^{T}_{y} (\boldsymbol{x}-\boldsymbol{y}) - \lVert \boldsymbol{x}-\boldsymbol{y} \rVert_{2}^{2},
\label{eq:quasiConvex}
\end{align}
for all $\boldsymbol{x},\boldsymbol{y}\in \mathbb{R}^{n}$, and $\boldsymbol{\zeta}_{y}\in \partial \lambda g(\boldsymbol{y})$.

The iterative procedure summarized in Algorithm \ref{alg:unrolling} is seen as 
\begin{align}
\boldsymbol{x}^{(t+1)} = \argmin_{\boldsymbol{x}\in \mathbb{R}^{n}} \hspace{0.5em} \left\langle \nabla f(\boldsymbol{x}^{(t)}),\boldsymbol{x}-\boldsymbol{x}^{(t)} \right\rangle + \frac{1}{2\alpha_{t} \lambda} \lVert \boldsymbol{x}-\boldsymbol{x}^{(t)} \rVert_{2}^{2} + g(\boldsymbol{x})
\label{eq:proofLem2}
\end{align}
We write equal in the above equation because the proximal in Eq. \eqref{eq:prox1} is invex therefore it always map to a global optimizer. From the Lipschitz continuous of $\nabla f$ we have
\begin{align}
f(\boldsymbol{x}^{(t+1)}) \leq f(\boldsymbol{x}^{(t)}) + \left\langle \nabla f(\boldsymbol{x}^{(t)}),\boldsymbol{x}^{(t+1)}-\boldsymbol{x}^{(t)} \right\rangle + \frac{L}{2} \lVert \boldsymbol{x}^{(t+1)}-\boldsymbol{x}^{(t)} \rVert_{2}^{2}.
\label{eq:proofLem3}
\end{align}
Considering the fact that from Eq. \eqref{eq:proofLem2} we conclude $-\nabla f(\boldsymbol{x}^{(t)}) + \frac{1}{\alpha_{t}\lambda} \left( \boldsymbol{x}^{(t)} -\boldsymbol{x}^{(t+1)}\right) \in \partial g(\boldsymbol{x}^{(t+1)})$, then Eq. \eqref{eq:quasiConvex} leads to 
\begin{align}
\lambda g(\boldsymbol{x}^{(t)}) - \lambda g(\boldsymbol{x}^{(t+1)}) &\geq \left\langle -\nabla f(\boldsymbol{x}^{(t)}) + \frac{1}{\alpha_{t}\lambda} \left( \boldsymbol{x}^{(t)} -\boldsymbol{x}^{(t+1)}\right) , \boldsymbol{x}^{(t)}-\boldsymbol{x}^{(t+1)} \right\rangle - \lVert \boldsymbol{x}^{(t)}-\boldsymbol{x}^{(t+1)} \rVert_{2}^{2} \nonumber\\
&\geq \left\langle \nabla f(\boldsymbol{x}^{(t)}) , \boldsymbol{x}^{(t+1)} - \boldsymbol{x}^{(t)} \right\rangle + \left(\frac{1}{\alpha_{t}\lambda} - 1\right)\lVert \boldsymbol{x}^{(t)}-\boldsymbol{x}^{(t+1)} \rVert_{2}^{2}.
\label{eq:proofLem4}
\end{align}
The above Eq. \eqref{eq:proofLem4} combined with Eq. \eqref{eq:proofLem3} yields
\begin{align}
\lambda g(\boldsymbol{x}^{(t)}) - \lambda g(\boldsymbol{x}^{(t+1)}) &\geq f(\boldsymbol{x}^{(t+1)}) - f(\boldsymbol{x}^{(t)}) - \frac{L}{2} \lVert \boldsymbol{x}^{(t+1)}-\boldsymbol{x}^{(t)} \rVert_{2}^{2} +\left(\frac{1}{\alpha_{t}\lambda} - 1\right)\lVert \boldsymbol{x}^{(t)}-\boldsymbol{x}^{(t+1)} \rVert_{2}^{2} \nonumber\\
f(\boldsymbol{x}^{(t)}) + \lambda g(\boldsymbol{x}^{(t)})&\geq f(\boldsymbol{x}^{(t+1)}) + \lambda g(\boldsymbol{x}^{(t+1)}) + \left( \frac{1}{\alpha_{t}\lambda} - 1 - \frac{L}{2} \right)\lVert \boldsymbol{x}^{(t)}-\boldsymbol{x}^{(t+1)} \rVert_{2}^{2}.
\label{eq:proofLem5}
\end{align}
Observe that by taking $\alpha_{t}<\frac{2}{L+2}$, then from Eq. \eqref{eq:proofLem5} we have that $f(\boldsymbol{x}^{(t)}) + \lambda g(\boldsymbol{x}^{(t)})\geq f(\boldsymbol{x}^{(t+1)}) + \lambda g(\boldsymbol{x}^{(t+1)})$ which is a sufficient decreasing condition. In addition, considering that the list of invex functions in Table \ref{tab:list}, and $f(\boldsymbol{x})=\lVert \boldsymbol{H}\boldsymbol{x} -\boldsymbol{v} \rVert_{2}^{2}$ are positive, then the loss function in Eq. \eqref{eq:problem4} is bounded below. Thus, in particular $f(\boldsymbol{x}^{(t)}) + \lambda g(\boldsymbol{x}^{(t)}) - (f(\boldsymbol{x}^{(t+1)}) + \lambda g(\boldsymbol{x}^{(t+1)})) \rightarrow 0$ as $t\rightarrow \infty$, which, combined with Eq. \eqref{eq:proofLem5}, implies that $\lVert \boldsymbol{x}^{(t)}-\boldsymbol{x}^{(t+1)} \rVert_{2}^{2}\rightarrow 0$ as $t\rightarrow \infty$. The later convergence implies the existence of fixed points to the proximal iteration in Algorithm \ref{alg:unrolling} (equivalently to Eq. \eqref{eq:proofLem2}), sufficient condition to guarantee that the sequence $\{\boldsymbol{x}^{(t)}\}$ convergences to a stationary point of $f(\boldsymbol{x}) + \lambda g(\boldsymbol{x})$. Thus, since in Theorem \ref{theo:ourCS} we proved the loss function in Eq. \eqref{eq:problem4} is invex then $\{\boldsymbol{x}^{(t)}\}$ converges to a global minimizer. 
\end{proof}
\newpage

\section{Image Compressive Sensing Experiments Evaluated with SSIM metric}
\label{app:newResults}
In this section we complement results of Experiments 1,2, and 3 of Section \ref{others}. We assess the imaging quality for these experiments using the structural similarity index measure (SSIM). The best and least efficient among invex functions is highlighted in boldface and underscore, respectively.

\begin{figure}[h]
	\centering
	\includegraphics[width=1\linewidth]{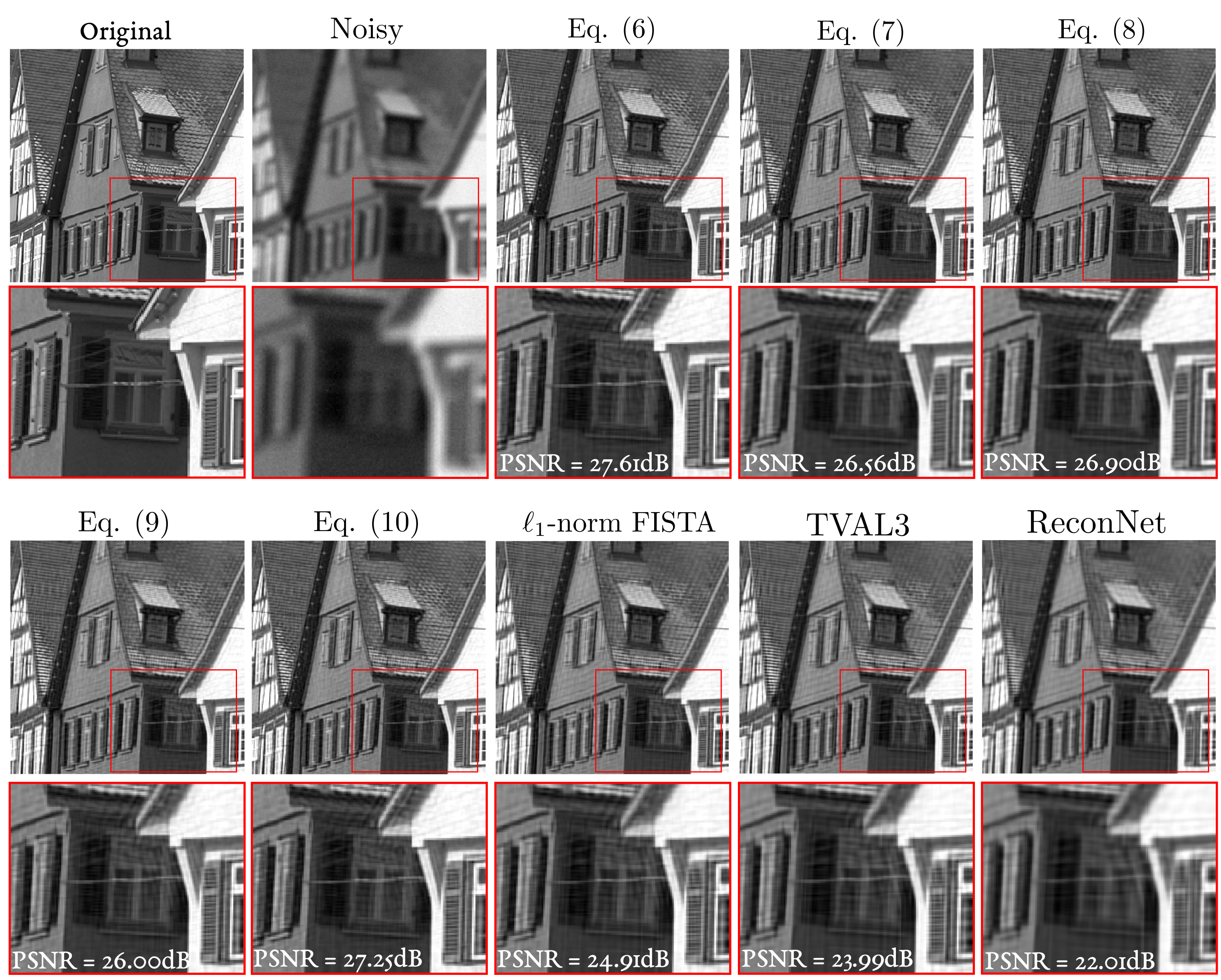}
	\vspace{-0.3em}
	\caption{Reconstructed images, for $SNR=30$dB, obtained by Algorithm \ref{alg:invexProximal} using Eqs. \eqref{fun1}, \eqref{fun2}, \eqref{fun3}, \eqref{fun4}, and \eqref{fun5}, which are compared with the outputs from FISTA, TVAL3, and ReconNet. In addition, to numerically evaluate their performance we estimate the PSNR for each image.}
	\label{fig:exp1}
\end{figure}

\textbf{Experiment 1} studies the effect of different invex regularizers under Algorithm \ref{alg:invexProximal}. The numerical results of this study are summarized in Table \ref{tab:globalResults1.1}. Also, we present Figure \ref{fig:exp1} which illustrates reconstructed images, for $SNR=30$dB, obtained by Eqs. \eqref{fun1}, \eqref{fun2}, \eqref{fun3}, \eqref{fun4}, and \eqref{fun5}, which are compared with the outputs from FISTA, TVAL3, and ReconNet. In addition, to numerically evaluate their performance we estimate the PSNR for each image.

\begin{table}[ht]
	\centering
	\caption{Comparison between convex and invex regularizers, in terms of SSIM, under Algorithm \ref{alg:invexProximal}, using $p=0.5$ for Eq.~\eqref{fun1}. }
	\resizebox{1\textwidth}{!}{\renewcommand{\arraystretch}{1.3}
		\begin{tabular}{P{1cm} | P{1.3cm} P{1.3cm} P{1.3cm} P{1.3cm} P{1.3cm} | P{2.5cm} P{2cm} P{2.2cm}}
			\hline
			& \multicolumn{5}{P{7.5cm}|}{(Experiment 1) Algorithm \ref{alg:invexProximal}, $p=0.5$ for Eq. \eqref{fun1}.} & FISTA \cite{beck2009fast} & TVAL3 \cite{li2013efficient} & ReconNet \cite{kulkarni2016reconnet}\\
			\hline
			SNR & Eq. \eqref{fun1} & Eq. \eqref{fun2} &Eq. \eqref{fun3} & Eq. \eqref{fun4} & Eq. \eqref{fun5} & $\ell_{1}$-norm & & \\
			\hline
			\centering $\infty$ & \textbf{0.9486} & 0.9370 & 0.9408 & $\underline{0.9332}$ & 0.9447 & 0.9257 & 0.9294 & 0.9220 \\
			\centering $20$dB & \textbf{0.8675} & 0.8495 & 0.8554 & $\underline{0.8437}$ & 0.8614 & 0.8323 & 0.8380 & 0.8267\\ 
			\centering $30$dB & \textbf{0.9055} & 0.8944 & 0.8981 & $\underline{0.8908}$ & 0.9018 & 0.8836 & 0.8872 & 0.8801\\
			\hline 			
		\end{tabular}
	}
	\label{tab:globalResults1.1}
\end{table}

\begin{figure}[h]
	\centering
	\includegraphics[width=1\linewidth]{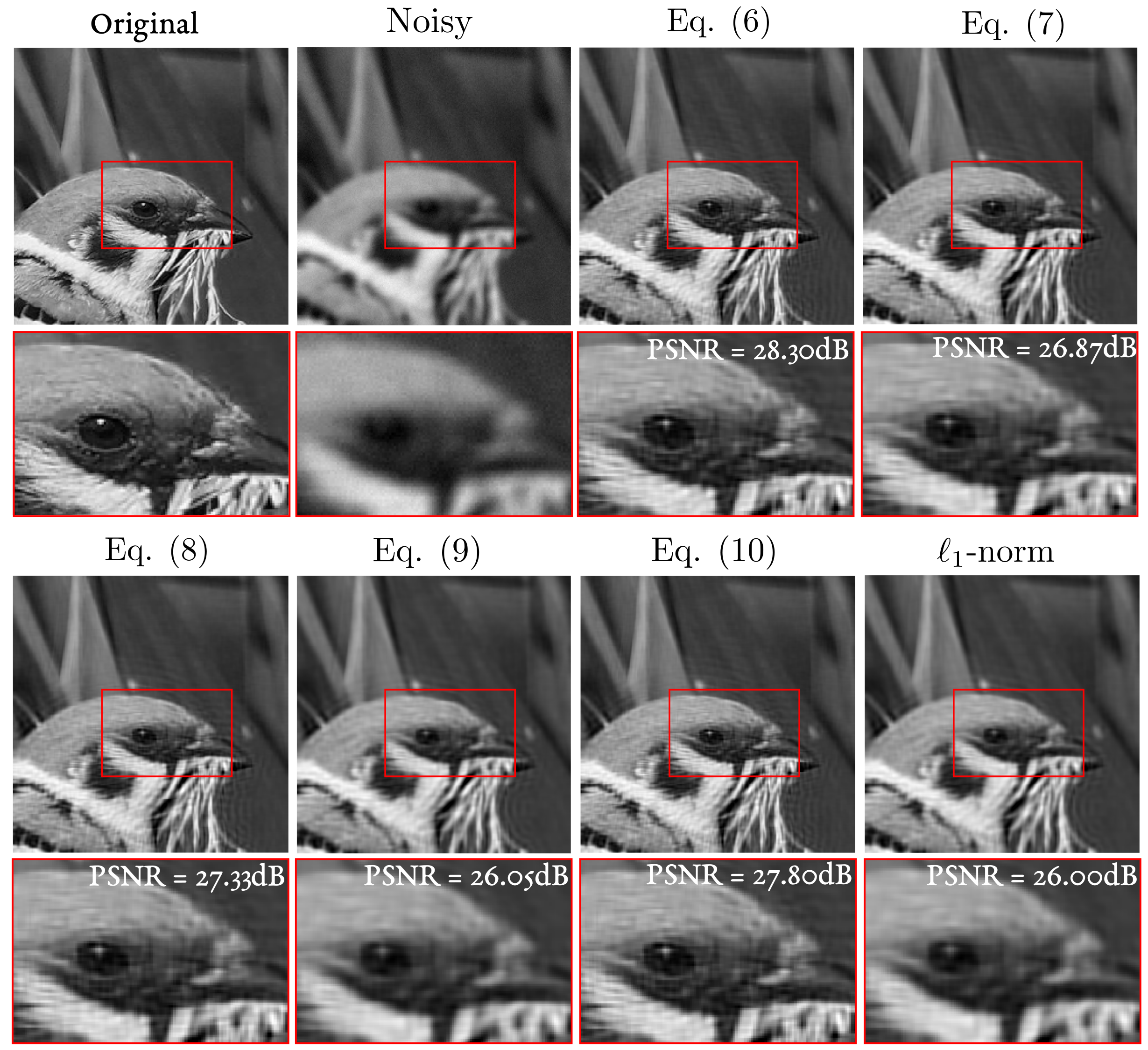}
	\vspace{-0.3em}
	\caption{Reconstructed images, for $SNR=30$dB, obtained by Algorithm \ref{alg:invexPnP} using Eqs. \eqref{fun1}, \eqref{fun2}, \eqref{fun3}, \eqref{fun4}, and \eqref{fun5}, which are compared with the outputs from $\ell_{1}$-norm. In addition, to numerically evaluate their performance we estimate the PSNR for each image.}
	\label{fig:exp2}
\end{figure}

\textbf{Experiment 2} studies the invex regularizers under the plug-and-play modification of Algorithm \ref{alg:invexProximal} as described in Section \ref{sub:PnP} \cite{krull2019noise2void}. The same deconvolution problem as in Experiment 1 is used. The numerical results of this study are summarized in Table \ref{tab:globalResults2.1}. Also, we present Figure \ref{fig:exp2} which illustrates reconstructed images obtained by Eqs. \eqref{fun1}, \eqref{fun2}, \eqref{fun3}, \eqref{fun4}, and \eqref{fun5}, which are compared with the outputs from $\ell_{1}$-norm. In addition, to numerically evaluate their performance we estimate the PSNR for each image.

\begin{table}[ht]
	\centering
	\caption{Comparison between convex and invex regularizers, in terms of SSIM, under plug-and-play Algorithm \ref{alg:invexPnP}, using $p=0.8$ for Eq. \eqref{fun1}.}
	\resizebox{1\textwidth}{!}{ \renewcommand{\arraystretch}{1.3}
		\begin{tabular}{P{1cm} | P{2cm} P{2cm} P{2cm} P{2cm} P{2cm} | P{2cm}}
			\hline
			& \multicolumn{5}{P{7.5cm}|}{(Experiment 2) Algorithm \ref{alg:invexPnP}, $p=0.8$ for Eq. \eqref{fun1}.} & \\
			\hline
			SNR &Eq. \eqref{fun1} & Eq. \eqref{fun2} & Eq. \eqref{fun3} & Eq. \eqref{fun4} & Eq. \eqref{fun5} & $\ell_{1}$-norm \\
			\hline
			\centering $\infty$ & \textbf{0.9581} & 0.9409 & 0.9465 & $\underline{0.9352}$ & 0.9523 & 0.9297\\
			\centering $20$dB & \textbf{0.8808} & 0.8680 & 0.8722 & $\underline{0.8638}$ & 0.8765 & 0.8597\\
			\centering $30$dB & \textbf{0.9189} & 0.9043 & 0.9091 & $\underline{0.8995}$ & 0.9140 & 0.8948\\
			\hline
		\end{tabular}
	}
	\label{tab:globalResults2.1}
\end{table}

\begin{figure}[H]
	\centering
	\includegraphics[width=0.95\linewidth]{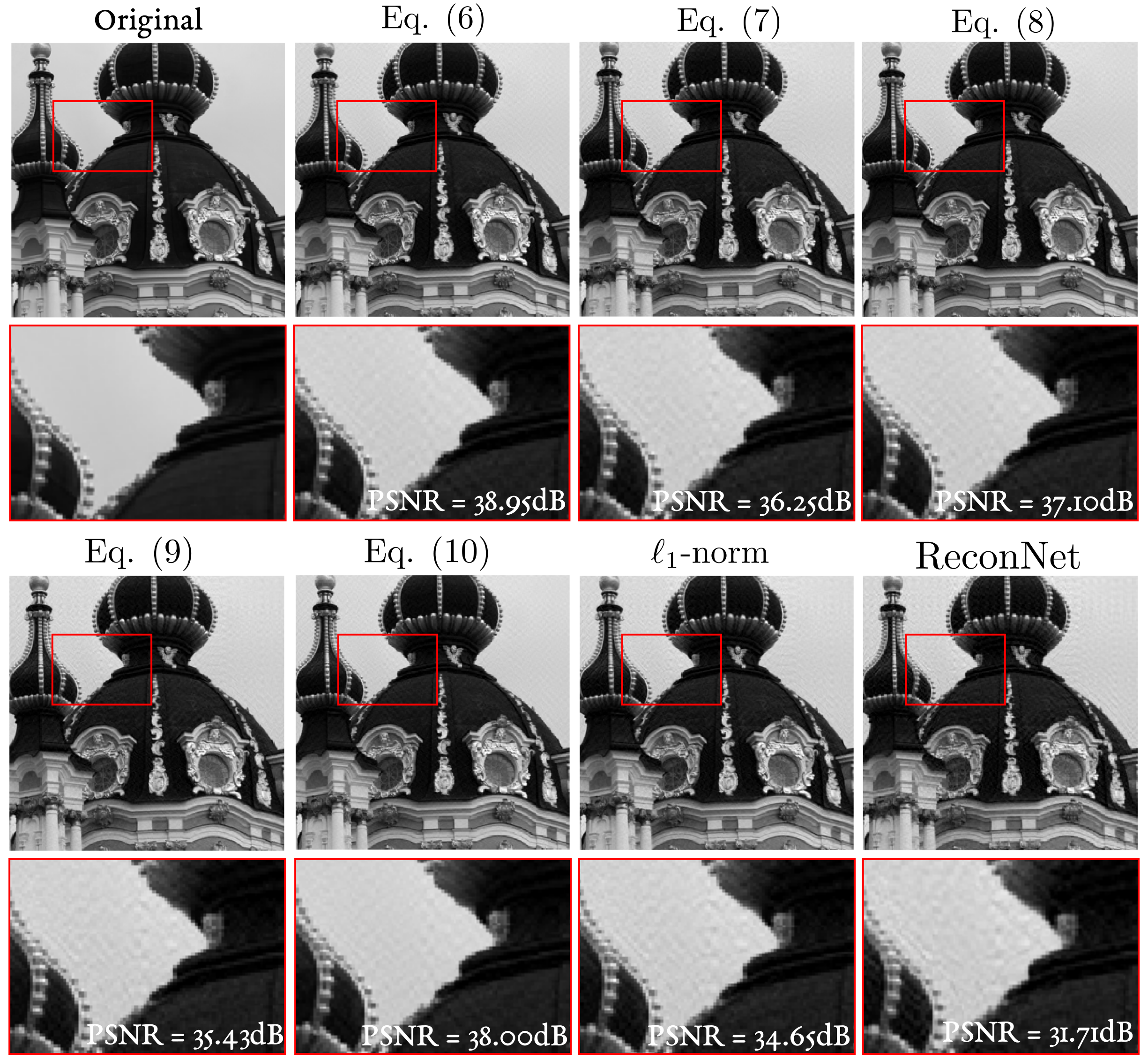}
	\vspace{-0.3em}
	\caption{Reconstructed images, for $SNR=30$dB, obtained by ISTA-Net using Eqs. \eqref{fun1}, \eqref{fun2}, \eqref{fun3}, \eqref{fun4}, and \eqref{fun5}, which are compared with the outputs from $\ell_{1}$-norm, where $m/n=0.6$. In addition, to numerically evaluate their performance we estimate the PSNR for each image.}
	\label{fig:exp4}
\end{figure}

\textbf{Experiment 3} compares the invex regularizers but under the unrolling framework as described in Section \ref{unrolling}. The numerical results of this study are summarized in Table \ref{tab:globalResults3.1}. Also, we present Figure \ref{fig:exp4} which illustrates reconstructed images obtained by Eqs. \eqref{fun1}, \eqref{fun2}, \eqref{fun3}, \eqref{fun4}, and \eqref{fun5} with ISTA-Net, which are compared with the outputs from $\ell_{1}$-norm + ISTA-Net. In addition, to numerically evaluate their performance we estimate the PSNR for each image.

\begin{table}[ht]
	\centering
	\caption{ Performance comparison between convex and invex regularizers, in terms of SSIM, for the unrolling experiment, using $p=0.85$ for Eq. \eqref{fun1}.}
	\resizebox{1\textwidth}{!}{ \renewcommand{\arraystretch}{1.3}
		\begin{tabular}{P{1cm}| P{1.1cm} | P{1.5cm} P{1.5cm} P{1.5cm} P{1.5cm} P{2.0cm} |P{2cm} | P{2.2cm}|}
			\hline
			& & \multicolumn{5}{P{10.5cm}|}{(Experiment 3) Algorithm \ref{alg:unrolling} - unfolded LISTA. $p=0.85$ for Eq. \eqref{fun1}} & LISTA \cite{chen2018theoretical}& ReconNet \cite{kulkarni2016reconnet}\\
			\hline
			SNR & $m/n$& Eq. \eqref{fun1} & Eq. \eqref{fun2} & Eq. \eqref{fun3} & Eq. \eqref{fun4} & Eq. \eqref{fun5} & $\ell_{1}$-norm & \\
			\hline
			\centering \multirow{3}{*}{$\infty$} & \centering 0.2 \vspace{0.1em} \hrule& \textbf{0.9279} & 0.9132 & 0.9181 & $\underline{0.9084}$ & 0.9230 & 0.9037 & 0.8990\\ 
			& \centering 0.4 \vspace{0.1em} \hrule& \textbf{0.9610} & 0.9423 & 0.9485 & $\underline{0.9363}$ & 0.9547 & 0.9303 & 0.9244\\
			& \centering 0.6 & \textbf{0.9890} & 0.9620 & 0.9708 & $\underline{0.9533}$ & 0.9798 & 0.9448 & 0.9364\\ 
			\cline{2-9}
			\centering \multirow{3}{*}{20dB} & \centering 0.2 \vspace{0.1em} \hrule& \textbf{0.8690} & 0.8628 & 0.8649 & $\underline{0.8608}$ & 0.8669 & 0.8587 & 0.8567 \\ 
			& \centering 0.4 \vspace{0.1em} \hrule& \textbf{0.9370} & 0.9205 & 0.9259 & $\underline{0.9151}$ & 0.9314 & 0.9098 & 0.9045\\ 
			& \centering 0.6 & \textbf{0.9498} & 0.9411 & 0.9440 & $\underline{0.9382}$ & 0.9469 & 0.9353 & 0.9325\\ 
			\cline{2-9}
			\centering \multirow{3}{*}{30dB} &\centering 0.2 \vspace{0.1em} \hrule& \textbf{0.8876} & 0.8781 & 0.8812 & $\underline{0.8750}$ & 0.8844 & 0.8719 & 0.8688 \\ 
			& \centering 0.4 \vspace{0.1em} \hrule& \textbf{0.9510} & 0.9318 & 0.9381 & $\underline{0.9255}$ & 0.9445 & 0.9194 & 0.9133 \\ 
			& \centering 0.6 & \textbf{0.9619} & 0.9545 & 0.9569 & $\underline{0.9520}$ & 0.9594 & 0.9496 & 0.9472\\ 
			\cline{1-9}
			& & \multicolumn{5}{P{10.5cm}}{(Experiment 3) Algorithm \ref{alg:unrolling} - unfolded ISTA-Net. $p=0.85$ for Eq. \eqref{fun1}} & & \\
			\cline{1-9}
			SNR & $m/n$& Eq. \eqref{fun1} & Eq. \eqref{fun2} & Eq. \eqref{fun3} & \multicolumn{1}{P{1.3cm}}{Eq. \eqref{fun4}} & \multicolumn{1}{P{1.3cm}|}{Eq. \eqref{fun5}}& $\ell_{1}$-norm \cite{zhang2018ista} & \\
			\cline{1-9}
			\centering \multirow{3}{*}{$\infty$} & \centering 0.2 \vspace{0.1em} \hrule& \textbf{0.9350} & 0.9219 & 0.9262 & $\underline{0.9176}$ & 0.9306 & 0.9134 & \multirow{9}{*}{-}\\ 
			& \centering 0.4 \vspace{0.1em} \hrule& \textbf{0.9733} & 0.9541 & 0.9604 & $\underline{0.9479}$ & 0.9668 & 0.9417 & \\
			& \centering 0.6 & \textbf{0.9899} & 0.9697 & 0.9763 & $\underline{0.9632}$ & 0.9831 & 0.9567 & \\ 
			\cline{2-8}
			\centering \multirow{3}{*}{20dB} & \centering 0.2 \vspace{0.1em} \hrule& \textbf{0.8829} & 0.8745 & 0.8773 & $\underline{0.8717}$ & 0.8801 & 0.8690 & \\ 
			& \centering 0.4 \vspace{0.1em} \hrule& \textbf{0.9501} & 0.9323 & 0.9382 & $\underline{0.9265}$ & 0.9441 & 0.9208 & \\ 
			& \centering 0.6 & \textbf{0.9611} & 0.9520 & 0.9550 & $\underline{0.9490}$ & 0.9580 & 0.9460 & \\ 
			\cline{2-8}
			\centering \multirow{3}{*}{30dB} &\centering 0.2 \vspace{0.1em} \hrule& \textbf{0.8990} & 0.8836 & 0.8887 & $\underline{0.8786}$ & 0.8938 & 0.8736 & \\ 
			& \centering 0.4 \vspace{0.1em} \hrule& \textbf{0.9641} & 0.9437 & 0.9504 & $\underline{0.9370}$ & 0.9572 & 0.9305 & \\ 
			& \centering 0.6 & \textbf{0.9859} & 0.9695 & 0.9749 & $\underline{0.9641}$ & 0.9804 & 0.9588 & \\ 
			\hline
		\end{tabular}
	}
	\label{tab:globalResults3.1}
\end{table}

\newpage

\section{Image Denoising Illustration}
\label{app:denoising}
For the sake of completeness we present in Algorithm \ref{alg:denoising} the denoising procedure employed in this paper (following \cite{cai2014data}) using invex regularizers $g(\boldsymbol{x})$ in Eqs. \eqref{fun1}, \eqref{fun3}, and \eqref{fun5}. The parameters $\lambda_{1},\lambda_{2}, K$, and $S$ were chosen to be the best for each analyzed function determined by cross validation.
\begin{figure}[h]
	\centering
	\includegraphics[width=1\linewidth]{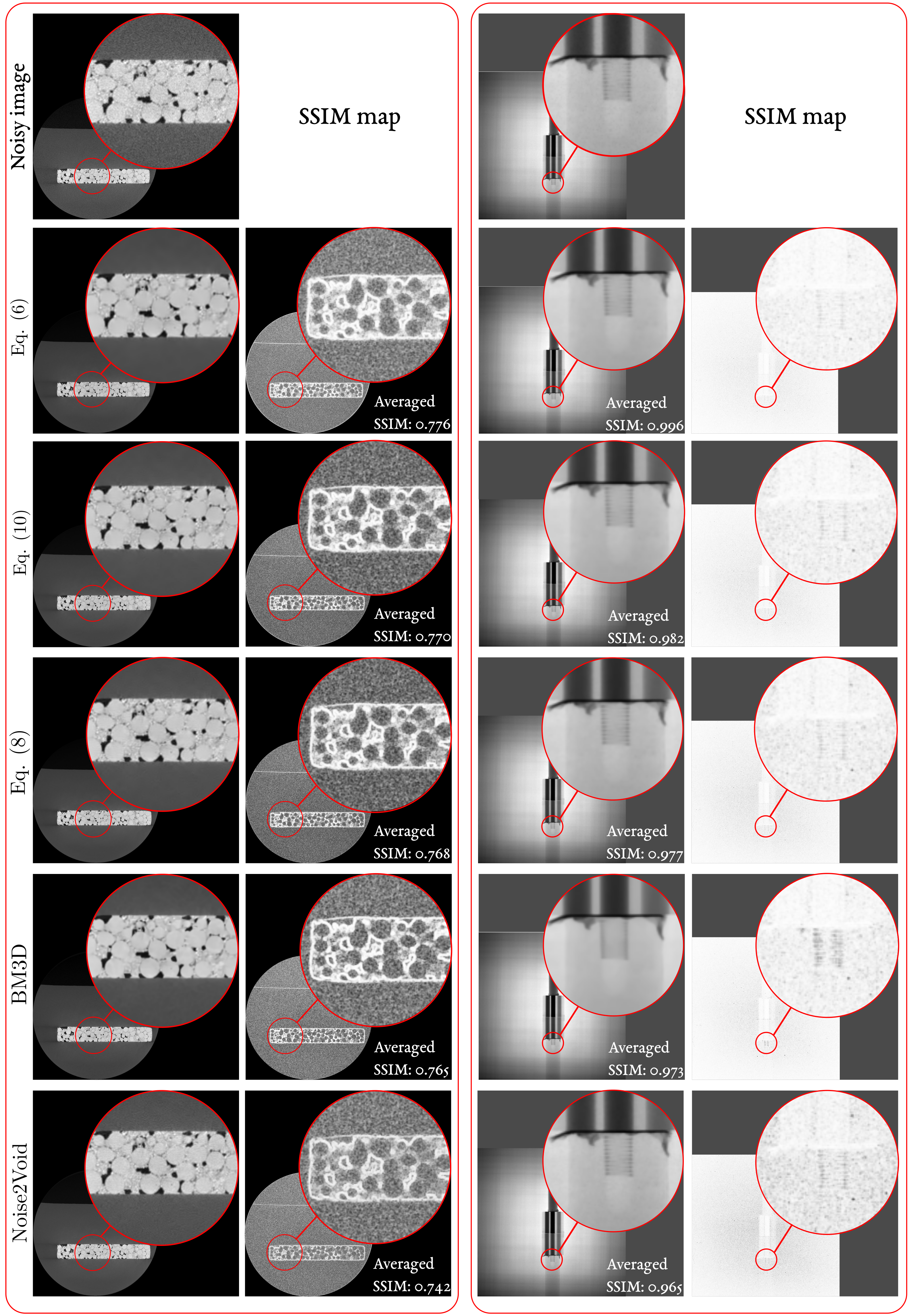}
	\vspace{-0.3em}
	\caption{Denoised image illustration for Eqs. \eqref{fun1},\eqref{fun3}, \eqref{fun5} and the state-of the-arts BM3D and Noise2Void. To evaluate the performance we employ the structural similarity index measure (SSIM) by reporting the SSIM map for each denoised image and its averaged value. Recall that SSIM is reported in the range $[0,1]$ where 1 is the best achievable quality and 0 the worst. In the SSIM map small values of SSIM appear as dark pixels. Thus, we conclude the best performance is achieved using the regularizer in Eq. \eqref{fun1} since it has the whitest SSIM maps.}
	\label{fig:denoising}
\end{figure}

\begin{algorithm}[ht]
	\caption{Denoising procedure using invex regularizers }
	\label{alg:denoising}
	\begin{algorithmic}[1]
		\State{\textbf{input}: noisy image $\boldsymbol{x}$, $S$ the number of patches of size $16\times 16$, $K$ number of iterations, and constant $\lambda_{1},\lambda_{2}\in (0,1]$.}
		\State{\textbf{initialize}: $\boldsymbol{W}^{(0)}=\frac{1}{256}\boldsymbol{1}$ where $\boldsymbol{1}\in \mathbb{R}^{256\times 256}$ is the matrix of ones.}
		\State{\textbf{Compute:} $\boldsymbol{P}\in \mathbb{R}^{256\times S}$ matrix containing random patches of size $16\times 16$ from $\boldsymbol{x}$}
		\State{$\boldsymbol{A} = (\boldsymbol{I}_{256}-\boldsymbol{W}^{(0)}(\boldsymbol{W}^{(0)})^{T})\boldsymbol{P}$, where $\boldsymbol{I}_{256}\in \mathbb{R}^{256\times 256}$ is the identity matrix}
		\For{$t=1$ to $K$}
		\State{$\boldsymbol{W}^{(t)}=(\boldsymbol{W}^{(t-1)})^{T}\boldsymbol{P}$}
		\State{$\hat{\boldsymbol{W}}^{(t)}[i,j]=\left \lbrace\begin{array}{ll}
			\boldsymbol{W}[i,j] & \lvert \boldsymbol{W}[i,j] \rvert \leq \lambda_{1}\\
			0 & \text{ otherwise }
			\end{array}\right.$}
		\State{run the SVD decomposition on $\boldsymbol{A}(\hat{\boldsymbol{W}}^{(t)})^{T}$ such that $\boldsymbol{A}(\hat{\boldsymbol{W}}^{(t)})^{T}=\boldsymbol{U}\boldsymbol{D}\boldsymbol{V}^{T}$.}
		\State{$\boldsymbol{W}^{(t)} = \boldsymbol{U}\boldsymbol{V}^{T}$}
		\EndFor
		\State{$\hat{\boldsymbol{x}}= (\boldsymbol{W}^{(K)})^{T}\text{Prox}_{\lambda_{2} g}(\boldsymbol{W}^{(K)}\boldsymbol{x})$\Comment{Denoising step}}
		
		\State{\textbf{return:} $\hat{\boldsymbol{x}}$\Comment{Denoised image}}
	\end{algorithmic}
\end{algorithm}

Employing Algorithm \ref{alg:denoising}, in Figure \ref{fig:denoising} we present some denoised images obtained by Eqs. \eqref{fun1}, \eqref{fun3}, \eqref{fun5}, which are compared with the outputs from BM3D, and Noise2Void. Since we are analyzing all the regularizers under non-ideal scenarios due to noise, results in Figure \ref{fig:denoising} highlight the benefit of having invex regularizers since the cleanest image is obtained by Eq. \eqref{fun1}. In addition, to numerically evaluate their performance we employ the structural similarity index measure (SSIM) by reporting the SSIM map for each denoised image and its averaged value. Recall that SSIM is reported in the range $[0,1]$ where 1 is the best achievable quality and 0 the worst. In the SSIM map small values of SSIM appear as dark pixels. Thus, we conclude the best performance is achieved using the regularizer in Equation \eqref{fun1} since it has the whitest SSIM maps (with highest SSIM values).
	
\end{document}